\documentclass[MSNbibl,number,citesort,seceqn,dvips]{arxbj}
\usepackage{mathbh,multirow,upgreek}
\usepackage{graphicx}

\aid{0}
\volume{19}
\issue{5B}
\pubyear{2013}
\firstpage{2652}
\lastpage{2688}
\doi{10.3150/12-BEJ470} %kopijuoti is PTS

\makeatletter
\newcommand{\rrvert}{\vert}
\newcommand{\llvert}{\vert}
\newcommand{\nref}[1]{(\ref{#1})}
\def\Cal{\mathcal}
\newcommand{\R}{\mathbb{R}}
\def\1{\mathbh{1}}
\newtheorem{theorem}{Theorem}[section]
\newtheorem{proposition}{Proposition}[section]
\newtheorem{lemma}{Lemma}[section]
\newremark{remark}{Remark}[section]

\newcommand{\CK}{{\Cal{K}}}
\newcommand{\CH}{{\Cal{H}}}
\newcommand{\CN}{{\Cal{N}}}
\newcommand{\CC}{{\Cal{C}}}
\newcommand{\CS}{{\Cal{S}}}
\newcommand{\eq}{\,{ \stackrel{\Delta}{=} }\,}
\def\a{\alpha}
\def\b{\beta}
\def\g{\gamma}
\def\Var{\operatorname{Var}}

\def\CK{\mathcal{K}}

\def\Var{\operatorname{Var}}
\makeatother

\begin{document}
\begin{frontmatter}

\title{Detection of a sparse submatrix of a high-dimensional noisy matrix}
\runtitle{Detection of a sparse submatrix}

\begin{aug}
%%%% inicialai - be tarpu
\author[1]{\fnms{Cristina} \snm{Butucea}\corref{}\thanksref{1}\ead[label=e1]{cristina.butucea@univ-mlv.fr}}%
\and
\author[2]{\fnms{Yuri I.} \snm{Ingster}\thanksref{2}\ead[label=e2]{yurii\_ingster@mail.ru}}%
\runauthor{C. Butucea and Y.I. Ingster} %% auto
\address[1]{Universit\'{e} Paris-Est, LAMA (UMR 8050), UPEMLV, UPEC, CNRS,
F-77464,
Marne-la-Vall\'ee, France and CREST, Timbre J340 3, av.
Pierre Larousse, 92240
Malakoff Cedex, France.\\
\printead{e1}}
\address[2]{St. Petersburg Electrotechnical University, 5,
Prof. Popov str., 197376 St. Petersburg, Russia.\\
\printead{e2}}
\end{aug}

% HISTORY:
\received{\smonth{9} \syear{2011}}
\revised{\smonth{7} \syear{2012}}

% ABSTRACT
%
\begin{abstract}
We observe a $N\times M$ matrix $Y_{ij}=s_{ij}+\xi_{ij}$ with
$\xi_{ij}\sim\CN(0,1)$ i.i.d. in $i,j$, and $s_{ij}\in\R$. We test
the null hypothesis $s_{ij}=0$ for all $i,j$ against the alternative
that there exists some submatrix of size $n\times m$ with
significant elements in the sense that $s_{ij}\ge a>0$. We propose
a test procedure and compute the asymptotical detection boundary $a$
so that the maximal testing risk tends to 0 as $M\to\infty$,
$N\to\infty$, $p=n/N\to0$, $q=m/M\to0$. We prove that this boundary
is asymptotically sharp minimax under some additional constraints.
Relations with other testing problems are discussed. We propose a
testing procedure which adapts to unknown $(n,m)$ within some given set
and compute the adaptive sharp rates. The
implementation of our test procedure on synthetic data shows
excellent behavior for sparse, not necessarily squared matrices. We
extend our sharp minimax results in different directions: first, to
Gaussian matrices with unknown variance, next, to matrices of
random variables having a distribution from an exponential family
(non-Gaussian) and, finally, to a two-sided alternative for matrices
with Gaussian elements.
\end{abstract}

% KEYWORDS
%
\begin{keyword}
\kwd{detection of sparse signal}
\kwd{minimax adaptive testing}
\kwd{minimax testing}
\kwd{random matrices}
\kwd{sharp detection bounds}
\end{keyword}

\end{frontmatter}
%
%%%%%%%%%%%%%%%%%%%%%%%%%%%
%s1 #&#
\section{Introduction}
%%%%%%%%%%%%%%%%%%%%%%%%%%

We observe a high-dimensional random matrix and we want to test
the occurrence of a particular submatrix of much smaller size,
which has elements with expected values larger than some threshold. We assume
that the entries of the matrix are independent, identically distributed (i.i.d.)
random variables but some underlying phenomenon can increase
significantly the
expected value of the random variables in the submatrix.

We have the observations that form an $N\times M$ matrix
$\mathbf{Y}=\{Y_{ij}\}_{i=1,\ldots,N, j=1,\ldots,M}$:
%
%e1.1 #&#
\begin{equation}
\label{mod} Y_{ij}=s_{ij}+ \sigma\xi_{ij},\qquad i=1,
\ldots,N, j=1,\ldots,M,
\end{equation}
where $\sigma>0$, $\{\xi_{ij}\}$ are i.i.d. random variables
and $s_{ij}\in\R$, for all $i\in\{1,\ldots,N\}$, $j\in\{1,\ldots,M\}$.
In the first part of the paper, the errors $\xi_{ij}$ are assumed to
have standard Gaussian law and $\sigma$ is assumed to be known. Without
loss of generality, we take $\sigma=1$ in this case. At the end of the
paper, we extend our results in different directions, as discussed
later on.
We test the null hypothesis that all elements of the matrix $\mathbf{Y}$
are i.i.d., standard Gaussian random variables $\CN(0,1)$, that is
%
%e1.2 #&#
\begin{equation}
\label{null} H_0\dvt \quad s_{ij}=0 \qquad\forall i=1,\ldots,N, j=1,
\ldots,M.
\end{equation}

The alternative under consideration will correspond to $n\times
m$-submatrices of sizes $n\in\{1,\ldots,N\}$, $m\in\{1,\ldots,M\}$ with
large enough entries. Let
%
%e1.3 #&#
\begin{equation}
\label{subset} A\subset\{1,\ldots,N\},\qquad \#(A)=n,\qquad B\subset\{1,\ldots,M\},\qquad \#(B)=m,\qquad
C=A\times B,
\end{equation}
and let ${\mathcal{C}}_{nm}$ be the collection of all subsets $C$ of the
form \nref{subset}. The set ${\mathcal{C}}_{nm}$ corresponds to the
collection of all $n\times m$ submatrices in $N\times M$ matrix. For
$a>0$, which may depend on $n,  m,   N$ and $M$.
We consider the alternative
%
%e1.4 #&#
\begin{equation}
\label{alt} H_1\dvt \quad \exists C\in{\mathcal{C}}_{nm} \mbox{ such
that } s_{ij}=0\mbox{ if } (i,j)\notin C \mbox{ and } s_{ij}
\ge a \mbox{ if } (i,j)\in C
\end{equation}
(in the Remark \ref{R1} below we discuss that a slightly larger
alternative can be considered).
The components of the matrix $\mathbf{Y}$ are independent under the
alternative as well. Denote by $P_S$ the probability measure that
corresponds to observations \nref{mod} with matrix $S=\{s_{ij}\}$
and by $E_S$ the expected value with respect to the measure $P_S$.

Let ${{\mathcal S}}_{nm,a}$ be the collection of all matrices $S=S_C$ that
satisfy \nref{alt}.

We discuss here only right-hand side alternatives, but, obviously,
left-hand side alternatives can be treated the same way for variables
$-Y_{ij}$ instead of $Y_{ij}$.

We extend our results to three different setups and sketch the proofs
of the results.
First, we consider errors having Gaussian distribution with unknown
variance $\sigma^2$.
We also consider other settings where the $Y_{ij}$'s come from an
exponential family.
Finally, in the initial case of Gaussian errors with known variance,
we consider a two-sided alternative of our test problem.

We are interested here in sparse matrices, that is, the case when $n$
is much
smaller than~$N$ and~$m$ is much smaller than $M$.

Sparsity assumptions were introduced for vectors. Estimation as well as
hypothesis testing
for vectors were thoroughly studied in the literature, see, for
example, Bickel, Ritov and Tsybakov~\cite{BRT}
and references therein and Donoho and Jin~\cite{DJ04}.

In the context of matrices, different sparsity assumptions can be imagined.
For example, matrix completion for low rank matrices with the nuclear
norm penalization
has been studied by Koltchinskii, Lounici and Tsybakov~\cite{KLT}.
Other results will be discussed later on.

We study the hypothesis testing problem under a minimax setting.
A test is any measurable function of the observations,
$\psi= \psi(\{Y_{ij}\})$ taking values in $[0,1]$.
For such a test $\psi=\psi(\{Y_{ij}\})$, we denote the probability of
type-I error, the probability of type-II error under simple alternative
and the maximal probability of type-II error over the set
${\mathcal S}_{nm,a}$ by
\[
\a(\psi)=E_{0}\psi,\qquad \b(\psi,S)=E_S(1-\psi),\qquad
\b_{nm,a}(\psi)=\sup_{S\in{{\mathcal S}}_{nm,a}}\b(\psi,S),
\]
respectively.
Let the risk be the following sum:
\[
\g(\psi,S)=\a(\psi)+\b(\psi,S),\qquad \g_{nm,a}(\psi)=\sup_{S\in
{{\mathcal S}}_{nm,a}}\g(
\psi,S)=\a(\psi)+\b_{nm,a}(\psi).
\]
We define the minimax risk at fixed level $\a\in(0,1)$ as
\[
\b_{nm,a, \a}=\inf_{\psi:\a(\psi)\le\a}\b_{nm,a}(\psi).
\]
Similarly, let the minimax testing risk be
\[
\g_{nm,a}=\inf_{\psi}\g_{nm,a}(\psi).
\]

From now on, we assume in the asymptotics that $N\to\infty, M\to\infty$ and $n=n_{NM}\to\infty,  m=m_{NM}\to\infty$. Other
assumptions will be given later.

We suppose that $a>0$ is unknown. The aim of this paper is to
give asymptotically sharp boundaries for minimax testing risk.
It means that, first, we are interested in the
conditions on $a=a_{NM}$ which guarantee distinguishability, that is,
the fact that $\g_{nm,a}\to0$
and $\b_{nm,a, \a}\to0$ for any $\a\in(0,1)$. We construct
a testing procedure based on a linear statistic combined with a scan statistic.
We prove the upper bounds of the minimax testing risk of this procedure.
Second, we describe conditions on $a$ for which we have
indistinguishability, that is, the convergence $\g_{nm,a}\to1$ and
$\b_{nm,a, \a}\to1-\a$
for any $\a\in(0,1)$. These results are called the lower bounds.
The two sets of conditions are complementary and match in rate and constant.

Often the sizes $n,m$ of submatrix are unknown, but we know a
set ${\mathcal{K}}_{NM}$ of couples of indices $(n,m) \in\{1,\ldots
,N\}\times
\{1,\ldots,M\}$ containing the true one.
Then we consider the ``adaptive'' problem for the combined
alternative ${{\mathcal S}}_{NM,{\mathbf{a}}}=\bigcup_{(n,m)\in
{\mathcal{K}}_{NM}}{{\mathcal S}}_{nm,a_{nm}}$, which corresponds to a
collection ${\mathbf{a}}=\{a_{nm}, (n,m)\in{\mathcal{K}}_{NM}\}$. The
quantities $\b_{NM,\mathbf{a},\a}, \g_{NM,{\mathbf{a}}}$ are defined in
a similar way as above. We define a testing procedure and check that,
if $a_{nm}$
satisfies the conditions for distinguishability uniformly over the collection
$\mathbf{a}$, the upper bounds still hold. The adaptive lower bounds
hold as
an easy consequence of the minimax lower bounds.

The problem of choosing a submatrix in a Gaussian random matrix has
been previously
studied by Sun and Nobel~\cite{SN}.
They were interested in maximal size submatrices of a matrix with
increasing size
in two setups. First, they consider the case when the average of the entries
of the submatrix is larger than
a given threshold and, second, when the entries are well-fitted by a two-way
ANOVA matrix in the least-squares sense
(i.e., the sum of squares of residuals is smaller than some given threshold).

The algorithm of choosing such submatrices was previously introduced in
Shabalin \textit{et al.}~\cite{SWPN},
who were also interested in finding large average submatrices.
This problem is strongly motivated by the research of gene expression
in microarray data.
In these large matrices, it is necessary to recover biclusters, that is
associations
between sets of samples (rows) and sets of variables (columns). These
associations
together with clinical and biological information are
``a first step in identifying disease subtypes and gene regulatory
networks''.
Many other algorithms for biclustering are discussed and compared on
real-data bases
concerning breast and lung cancer studies.

Similar problems were considered in Addario-Berry \textit{et al.} \cite{comb}.
They use the same testing procedures for
vectors of random variables, where the alternatives may have various
combinatorial structures.
In particular, they consider the example of detecting a clique of a
certain size
in a graph and they compute upper and lower bounds for the Bayesian
test error.
A bipartite graph of size $(N,M)$ is a graph having edges only between
the $N$
vertices of one set to the $M$ vertices of a second set.
A biclique is a complete bipartite subgraph of size $(n,m)$, that is,
a subgraph where
all $n$ vertices from the first set are connected to the $m$ vertices
from the second set.
We consider the problem of detecting a biclique. Our results are sharp minimax
and adaptive to the size of the unknown biclique.\looseness=1

The plan of the paper is as follows. In Section~\ref{cons}, we give the
test procedures.
We state the conditions on the detection boundary $a$ such that
distinguishability
is possible.
Under mild additional assumptions, we give the conditions on $a$ so
that the alternative is
indistinguishable from the null hypothesis.

In Section~\ref{adap}, we consider the adaptive setup where $(n,m)$ is
unknown but
belongs to some collection of sequences $\mathcal{K}_{NM}$.
We compute the adaptive rates of testing
of a slightly modified test procedure.

In Section~\ref{simu}, we perform a numerical study of the procedures
that attain
the sharp upper bounds.
In order to compute the scan statistic, a heuristic stochastic
algorithm from Shabalin
\textit{et al.}~\cite{SWPN} is used. The empirical detection boundary is
very close to the
one predicted by our results.
%For a fixed first-type error of $1\%$, the change of regime in the
%second-type error occurs very tightly around the detection boundary
%(i.e. it is very close to 1
%before the boundary and decreases abruptly to 0 after the detection
%boundary).

In Section~\ref{ext}, we give extensions of our results to Gaussian
variables of unknown variance $\sigma^2$, to non-Gaussian matrices
with distribution in an exponential family and to
two-sided tests for Gaussian matrices, respectively.

We include in Section~\ref{resadd} comments to understand how our
results compare to previously studied alternatives:
subsets without structure and rectangular submatrices.
The first case can be assimilated to detection of a sparse signal in vector
observations of length $N\times M$, so the set of alternatives and
the detection boundary are much larger than in our case.
We summarize well-known results by Ingster~\cite{I97},
Ingster and Suslina~\cite{IS02b} and Donoho and Jin~\cite{DJ04}.
The second case is the detection of rectangles in the large matrix
(connected submatrices), which constitutes a set of alternatives
smaller than ours.
This case is studied in Arias-Castro \textit{et al.}~\cite{ACDH} and
\cite
{ACCD11} for other geometric shapes of clusters. In order to be self-contained,
we state and prove sharp upper and lower bounds, for the rectangular clusters.

Section~\ref{PL} is mainly concerned with the proof of the lower bounds
stated in Section~\ref{ML}.
The \hyperref[app]{Appendix} contains the proofs of the other results of the paper.

%%%%%%%%%%%%%%%%%%%%%%%%%%%%%%%%%%%
%s2 #&#
\section{Main results}\label{M}
%%%%%%%%%%%%%%%%%%%%%%%%%%%%%%%%%%%

We denote by $n=n_{NM},   m=m_{NM}$ and $a=a_{N,M}$.

Denote also $p=n/N, q=m/M$. From now on, we suppose that
%
%e2.1 #&#
\begin{equation}
\label{assump1} N \to\infty,\qquad M \to\infty,\qquad n \to\infty,\qquad m \to\infty\qquad \mbox{such
that } p\to0, q\to0.
\end{equation}

For general sequences $\{u_n\}_{n \geq1}$ and $\{v_n\}_{n \geq1}$
of real numbers, such that $v_n >0$ for $n$ large enough, we say
that the sequences are asymptotically equivalent, $u_n \sim v_n$, if
$\lim_{n \to\infty} u_n/v_n = 1$. Moreover, we say that the
sequences are asymptotically of the same order, $u_n \asymp v_n$, if
there exists two constants $0<c \leq C <\infty$ such that $c\leq
\liminf_{n \to\infty} u_n/v_n $ and $\limsup_{n \to\infty} u_n/v_n
\leq C$.

%s2.1 #&#
\subsection{Known size of the submatrix}\label{cons}

In a minimax setup, we suppose that for each $N,  M$ we know $n$ and $m$.

Let us consider two test procedures, one based on a linear statistic
$\psi^{\mathrm{lin}}_H$ and
the other based on a scan statistic $\psi^{\mathrm{max}}$. The final test
procedure $\psi^*$
will reject as soon as at least one of them rejects the null hypothesis.

%s2.1.1 #&#
\subsubsection{Test procedure and its performance}

The first test procedure $\psi^{\mathrm{lin}}_H$ is based on the linear statistic
\[
t_{\mathrm{lin}}=\frac{1}{\sqrt{NM}}\sum_{i,j}Y_{ij},\qquad
\psi^{\mathrm{lin}}_H=\1_{t_{\mathrm{lin}}>H}.
\]
%
%One easily gets the following non-asymptotic result.
%
The second test $\psi^{\mathrm{max}}$ is based on the maximal sum over all
submatrices. Put
%
%e2.2 #&#
\begin{equation}
\label{Y} Y_C=\frac{1}{\sqrt{nm}}\sum_{(i,j)\in C}Y_{ij},
\end{equation}
and
%
%e2.3 #&#
\begin{equation}
\label{psimax} t_{\mathrm{max}}=\max_{C\in{{\mathcal C}}_{nm}}Y_C,\qquad
\psi^{\mathrm{max}}=\1_{t_{\mathrm{max}}>T_{nm}},
\end{equation}
where $T_{nm}=\sqrt{2\log(G_{nm})}, G_{nm}=\#({\mathcal C}_{nm})
= {N\choose n} {M\choose m}$.
The computation of this statistic is discussed in Section~\ref{simu}.

The following theorem gives sufficient conditions for the detection
boundary $a$
such that distinguishability holds. The test procedure which attains
these bounds is
\[
\psi^* = \max\bigl\{\psi^{\mathrm{lin}}_H, \psi^{\mathrm{max}}\bigr
\}
\]
for properly chosen $H$.

%th2.1 #&#
\begin{theorem}[(Upper bounds)] \label{TU}
Assume \nref{assump1} and let $a$ be such that at least one of
the following conditions hold
%
%e2.4 #&#
\begin{equation}
\label{cond1} a^2nmpq\to\infty
\end{equation}
or
%
%e2.5 #&#
\begin{equation}
\label{cond2} \liminf\frac{a^2nm}{2(n\log(p^{-1})+m\log(q^{-1}))}>1.
\end{equation}
Then $\psi^*$ with $H\to\infty$ and such that $H \leq ca\sqrt {nmpq}, c<1$ when \nref{cond1} holds,
%and with $T_{nm} = \sqrt{2 \log G_{nm}}$, $G_{nm} = {N \choose n}
satisfies $\gamma_{nm,a}(\psi^*) \to0$.
\end{theorem}

Proof is given in Appendix~\ref{PPmax}.
%The test procedure $\psi^*$ rejects the null hypothesis as soon as
%either the linear or the scan test rejects.

Formally, the procedure has a simple structure. Nevertheless, there are
difficulties for computation of the scan statistic in the matrix case.
Indeed, in the vector case, it is enough to order increasingly all the
elements and take the sum of the largest values. In the matrix case, we
have no such simple ordering.
We shall discuss in the numerical study below the empirical algorithm
used to compute the scan statistic.

Let us also note that this procedure assumes that $n$ and $m$ are
known. An adaptive version of the scan test
will be given in the next section.

%s2.1.2 #&#
\subsubsection{Lower bounds}\label{ML}

In this section, we obtain matching lower bounds that apply to all
tests under additional assumptions on the matrix and submatrix sizes.
We discuss these assumptions after the theorem.
%
%th2.2 #&#
\begin{theorem}[(Lower bounds)]\label{TL}
Assume \nref{assump1} and
%
%e2.6 #&#
\begin{equation}
\label{cond3a} \frac{\log\log(p^{-1})}{\log(q^{-1})}\to0, \qquad\frac{\log\log(q^{-1})}{\log(p^{-1})}\to0.%,
\end{equation}
Moreover, assume that
%
%e2.7 #&#
\begin{equation}
\label{aa} {n}\log\bigl(p^{-1}\bigr)\asymp m\log\bigl(q^{-1}
\bigr),
\end{equation}
and that the following two conditions are satisfied:
%
%e2.8 #&#
\begin{equation}
\label{cond12oppa} a^2nmpq\to0
\end{equation}
and
%
%e2.9 #&#
\begin{equation}
\label{cond12oppb} \limsup\frac{a^2nm}{2( n\log(p^{-1})+m\log(q^{-1}))}<1.
\end{equation}
Then the distinguishability is impossible, that is, $\g_{nm,a}\to1 $ and
$\b_{nm,a,\a}\to1-\a$ for any $\a\in(0,1)$.
\end{theorem}

Proof is given in Section~\ref{PL}.

These results for the upper and the lower bounds can be interpreted as
follows. Under the conditions \nref{assump1}, \nref{cond3a} and \nref{aa},
a sharp detection boundary $a^*$ is defined via the relations
%
%e2.10 #&#
\begin{equation}
\label{cond4} \bigl(a^*\bigr)^2nmpq\asymp1,\qquad \bigl(a^*
\bigr)^2nm\sim2\bigl(n\log\bigl(p^{-1}\bigr)+m\log
\bigl(q^{-1}\bigr)\bigr),
\end{equation}
in the problem with known $(n,m)$.
Note that the detection boundary can be written as
\[
a^* = \min \biggl\{ \frac{1}{\sqrt{nmpq}}, \sqrt{\frac{2(n\log(p^{-1})+m\log(q^{-1}))}{nm}} \biggr\} .
\]

The additional assumptions \nref{cond3a} and \nref{aa} appearing in
the previous lower bounds are satisfied, for example, in the case
where $n\sim c m$, for some $0< c < \infty$, and for $N\sim n^A$ and
$M\sim m^B$ for $A$ and $B$ larger than 1. In this case, the
detection boundary is of the form:
\begin{eqnarray*}
a^*&\asymp& n^{-2+(A+B)/2}\qquad\mbox{if } A+B\le3,
\\
a^*&\sim&\sqrt{\frac{2D\log(n)}{n}}\qquad \mbox{if } A+B> 3, \mbox{where } D=(A-1)c+B-1.
\end{eqnarray*}
The particular case when $A=B>1, c=1$ is the case of asymptotically
squared matrices and submatrices, and we get
\begin{eqnarray*}
a^*&\asymp& n^{-2+A}\qquad\mbox{if } A\le3/2,
\\
a^*&\sim& 2\sqrt{\frac{(A-1)\log(n)}{n}}\qquad\mbox{if } A> 3/2.
\end{eqnarray*}

%re2.1 #&#
\begin{remark}\label{R1}
We can state the alternative hypothesis in a more general form:
\[
H_1\dvt\quad  \exists C \in\mathcal{C}_{nm} \mbox{ such that }s_{ij} = 0 \mbox { if } (i,j) \notin C \mbox{ and } \sum
_{(i,j) \in C} s_{ij} \geq a nm.
\]
Indeed, our probabilities of error depend on the elements of the
submatrix $C$ only through the sum of its elements. Therefore, the previous
test procedure will attain the same rates and the same lower bound techniques
will give the previous results for this more general test problem.
\end{remark}

%s2.2 #&#
\subsection{Adaptation to the size of the submatrix}\label{adap}

If the size $(n,m)$ of the submatrix $C$ with significantly large
elements under
the alternative (\ref{alt}) is unknown, we suppose that
it belongs to the set $\mathcal{K}_{NM}$, for each $N$ and $M$.
The alternative hypothesis can be written
\begin{eqnarray*}
H_1(\mathcal{K}_{NM})\dvt&&  \quad \exists(n,m)\in\mathcal {K}_{NM},
\exists C\in{\mathcal{C}}_{nm} \mbox{ such that}
\\
&& \quad s_{ij}=0, \mbox{ if } (i,j)\notin C \mbox{ and } s_{ij}\ge
a_{nm}, \mbox{ if } (i,j)\in C.
\end{eqnarray*}
Additionally, we suppose that the sequence of sets $\{\mathcal
{K}_{NM}\}_{N,M}$ is such that
\[
\sup_{(n,m)\in{\mathcal{K}}_{NM}} \biggl( \frac1n + \frac1m +\frac nN +\frac mM \biggr) \to0
\]
as $N,  M \to\infty$.

This implies that
%
%e2.11 #&#
\begin{equation}
\label{assump2} %\tau({\mathcal{K}}_{NM})\eq
\sup_{(n,m)\in{\mathcal{K}}_{NM}} \biggl(\frac{\log(N)}{n\log(p^{-1})}+
\frac{\log(M)}{m\log(q^{-1})} \biggr) %\tau(N,M,n,m)
\to0\qquad \mbox{as } N, M \to\infty.
\end{equation}
The set $\mathcal{K}_{NM}$ contains sizes of submatrices that we have
to explore
in order to test in an adaptive way. Therefore, previous assumption
insure, on the one hand, that $p \to0$ and $q\to0$ uniformly over
$(n,m) \in\mathcal{K}_{NM}$ as $N, M \to\infty$ and, on the other
hand, that the least size of the submatrices still grows to infinity
with $N$ and $M$.

The adaptive test procedure is $\psi_{NM}^* = \max\{\psi^{\mathrm{lin}}_H,
\psi^{\mathrm{max}}_{NM}\}$,
where $\psi_H^{\mathrm{lin}}$ is the linear statistic defined in Section~\ref{cons}
and $\psi_{NM}^{\mathrm{max}}$ is a modified version of $\psi^{\mathrm{max}}$ defined
as follows.
Indeed, the linear statistic is free of $n$ and $m$, but the scan
statistic is not and, therefore, normalization will occur for each
possible $(n,m)$.
Set
\[
V_{nm}=\sqrt{2\log(NMG_{nm})},\qquad\!\! t_{NM,\mathrm{max}}=
\max_{(n,m)\in\CK_{NM}}\max_{C\in{{\mathcal C}}_{nm}}Y_C/V_{nm},\qquad\!\!
\psi^{\mathrm{max}}_{NM}=\1_{t_{NM,\mathrm{max}}>1}.
\]

The adaptive test will reject the null hypothesis as soon as at least
one between the linear test or the scan tests associated to each $(n,m)
\in\mathcal{K}_{NM}$ rejects.

%th2.3 #&#
\begin{theorem}\label{TLa} Assume \nref{assump1} and let the set
${\mathcal{K}}_{NM}$
be such that condition \nref{assump2} holds.

\emph{Upper bounds}. Let $\mathbf{a}=\mathbf{a}_{NM}=\{a_{nm},
(n,m)\in{\mathcal{K}}_{NM}\}$ be detection boundaries such that at
least one of the following conditions hold
%
%e2.12 #&#
\begin{equation}
\label{cond1a} \min_{(n,m)\in{\mathcal{K}}_{NM}}a_{nm}^2nmpq \to\infty
\end{equation}
or
%
%e2.13 #&#
\begin{equation}
\label{cond2a} \liminf\min_{(n,m)\in{\mathcal{K}}_{NM}} \frac{a_{nm}^2nm}{2(n\log(p^{-1})+m\log(q^{-1}))}>1.
\end{equation}

Then, $\psi^*_{NM}$ with $H\to\infty$ such that $H \leq c \min_{(n,m)\in{\mathcal{K}}_{NM}}a_{nm} \sqrt{nmpq}$ for some $0<c<1$ when
\nref{cond1a} holds, is such that $\gamma_{NM,\mathbf{a}}(\psi^*_{NM})
\to0$.
\end{theorem}

Proof is given in Appendix~\ref{proofTLa}.

The previous theorem actually shows that the test procedure is adaptive
to the size $(n,m)$ of the submatrix as far as the assumptions hold
uniformly. Indeed, the linear procedure is free of the size of the
submatrix and the scan statistic adapts to $(n,m)$ without any loss in
the rate.

The lower bounds in the adaptive setup are an obvious consequence of
Theorem~\ref{TL}. Let us state the adaptive lower bounds:
Suppose that for each $N$, $M$ there exists
$(n^*,m^*)$ in the collection ${\mathcal{K}}_{NM}$ such that
\[
\frac{\log\log(N/n^*)}{\log(M/m^*)}\to0,\qquad
\frac{\log\log(M/m^*)}{\log(N/n^*)}\to0
\]
and that
$
{n^*}\log(N/n^*)\asymp m^*\log(M / m^*),
$
as $N\to\infty$ and $M\to\infty$.
Let $\mathbf{a}=\mathbf{a}_{NM}= \{a_{nm}, (n,m)\in{\mathcal{K}}_{NM}\}$ be
such that
\[
a_{n^*m^*}^2n^*m^*p^*q^* \to0
\]
and
\[
\limsup %\min_{(n,m) \in{\mathcal{K}}_{NM}}
\frac{a_{n^*m^*}^2n^*m^*}{2 (n^*\log({p^*}^{-1})+m^*\log({q^*}^{-1}))}<1.
\]
Then $\g_{NM,\mathbf{a}}\to1 $ and $\b_{NM,\mathbf{a},\a}\to1-\a$ for any
$\a\in(0,1)$.

%%%%%%%%%%%%%%%%%%%%%%%%
%s3 #&#
\section{Simulations} \label{simu}
%%%%%%%%%%%%%%%%%%%%%%%%

We have implemented the testing procedure $\psi^* = \max\{\psi^{\mathrm{lin}}_H,\psi^{\mathrm{max}}\}$
on synthetic data.
While the linear procedure is rather obvious, the computation of the
statistic $t_{\mathrm{max}}=\max_{C \in\CC_{nm}} Y_C$ is done by using the heuristic
algorithm introduced and studied empirically by Shabalin \textit{et al.}
\cite{SWPN}.
This algorithm is also implemented and studied by Sun and Nobel~\cite
{SN} with good empirical results.

Let us briefly recall this algorithm: we choose randomly a
set of $n$ rows out of $N$. Then, we sum in every column
the elements of the previously selected rows. We select now the
columns corresponding to the $m$ largest sums obtained in this way.
We sum, next, in every row the elements belonging to the selected
columns and
select the rows corresponding to the $n$ largest sums.
We repeat the algorithm until the sum of elements $Y_{ij}$ of the
selected submatrix
does not increase anymore.
As the procedure can stop at a local maximum, we repeat the procedure
$K$ times,
where $K$ is large (in our simulation $K=10\,000$).
We take the maximum value of the outputs.
This replication is needed to enforce that with high probability the
output approaches
the global maximum.

%f1 #&#
\begin{figure}

\includegraphics{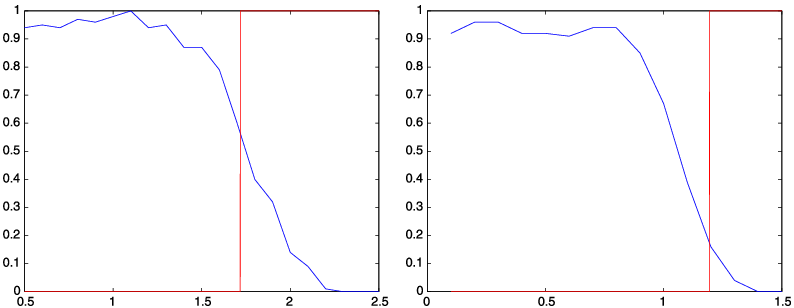}

\caption{Estimated second-type error probability for fixed $\a=1\%$,
detection boundary $a^*$, $N=M=200$; for $n=m=5$, $a^*= 1.7179$ (left),
for $n=m=10$, $a^*= 1.1943$ (right).}\label{figure1}
\end{figure}

%f2 #&#
\begin{figure}

\includegraphics{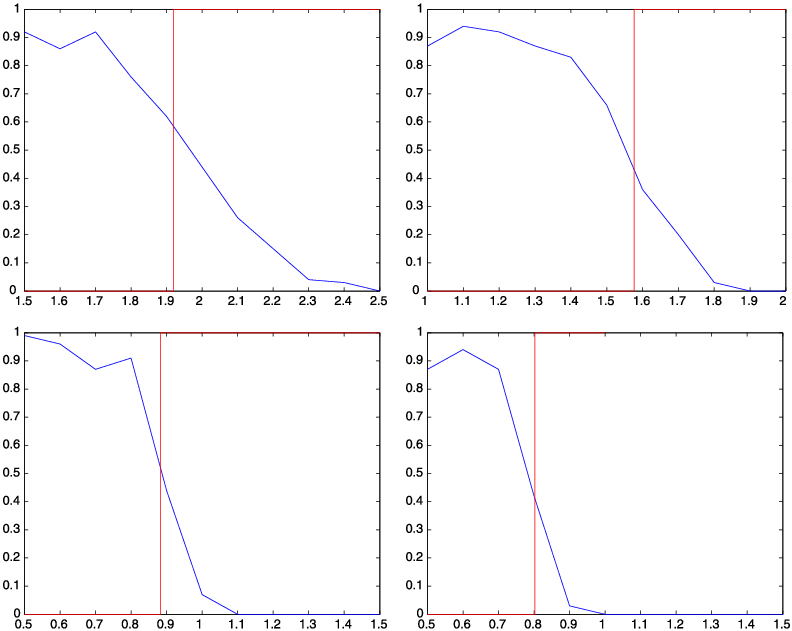}

\caption{Estimated second-type error probability for fixed $\a=1\%$,
detection boundary $a^*$, $N=M=500$; for $n=m=5$, $a^*= 1.9194$ (upper left),
for $n=5$, $m=10$, $a^*= 1.5767$ (upper right), for $n=15$, $m=20$,
$a^* = 0.8831$ (lower left), for $n=m=20$, $a^* = 0.8024$ (lower
right).}\label{figure2}
\end{figure}

We have simulated matrices of size $ N \times M$ of i.i.d. standard Gaussian
random variables for $N=M=200$ and $N=M=500$.

We calibrated the test statistics $\psi_H^{\mathrm{lin}}$ and $\psi^{\mathrm{max}}$ in such
a way that the type-I error occurs with probability $\alpha(\psi^*)
\leq1\%$.
This calibration is done by using the Gaussian
quantile $H=2.3262$ for $\psi^{\mathrm{lin}}_H$ and the empirical quantile (out
of 100 samples)
for $\psi^{\mathrm{max}}$.

Then, we have added the value $a>0$ to the elements of the upper
left submatrix of size $n\times m$. From resulting observations,
we compute $\psi^* = \max\{\psi_H^{\mathrm{lin}},\psi^{\mathrm{max}}\}$. We repeat the
test $L=100$ times
and average the values of the test procedure $\psi^*$. Denote by $\bar
\psi^*$ this average
and note that $1-\bar\psi^*$ estimates the probability of type-II error.

We plot the estimated second-type error probabilities for different
values of $a$ in the neighborhood of the detection boundary predicted
by our theorems, for different values of $n$ and $m$. The results in
Figure~\ref{figure1} correspond to $N=M=200$, while in
Figure~\ref{figure2} to $N=M=500$.

Figures \ref{figure1} and \ref{figure2} show that the empirical
detection boundary is very close to $a^*$ which is predicted by out
theoretical results. Indeed, the second-type error probability is close
to 0.5 at some point close to $a^*$. The plots also show very fast
decay of this probability on a small vicinity of $a^*$. This means that
the test is very powerful for values of $a$ slightly larger than the
detection boundary $a^*$. Note also that, for fixed $N$ and $M$, $a^*$
decreases to 0 as $n$ and $m$ increase.

%%%%%%%%%%%%%%%%%%%%%%%%%%%%%%%%%%%%%%%%%%%%%%%%%%%%
%s4 #&#
\section{Extensions}\label{ext}
%%%%%%%%%%%%%%%%%%%%%%%%%%%%%%%%%%%%%%%%%%%%%%%%%%%%

We extend our results in different directions. First, we consider
matrices of
i.i.d. random variables having Gaussian law with unknown variance
$\sigma$,
next, random variables having a distribution belonging to
the exponential family (not necessarily Gaussian)
and, finally, test problem with two-sided alternative for the Gaussian matrices.

%s4.1 #&#
\subsection{Extension to Gaussian variables with unknown variance}

Sharp results in Theorems~\ref{TU} and \ref{TL} still hold if the random
variables $Y_{ij}$ have unknown variance $\sigma$, under a mild additional
assumption. We sketch here the test procedure and proof of the upper bounds.

We estimate the unknown variance $\sigma^2$ of our data by $\hat
\sigma^2$, where
\[
\hat\sigma^2 = \frac1{NM} \sum_{i,j}
Y_{ij}^2.
\]
This estimator is unbiased under the null hypothesis, but biased under
the alternative.\vadjust{\goodbreak}

We replace $Y_{ij}$ by $Y_{ij}/\hat\sigma$ and slightly enlarge
$T_{nm}$ in the test procedure $\psi^*$. We denote by $\hat t_{\mathrm{lin}}
= t_{\mathrm{lin}}/\hat\sigma$, $\hat t_{\mathrm{max}} = t_{\mathrm{max}}/\hat\sigma$ and put
\[
\hat\psi^* = \max\{ \1_{ \hat t_{\mathrm{lin}}>H }, \1_{\hat
t_{\mathrm{max}}>T_{nm,\delta}}\}
\]
for $T_{nm,\delta}=\sqrt{ (2+\delta) \log\bigl({N \choose n} {M \choose
m}\bigr)}$ and some $\delta>0$ small enough. Recall that $T_{nm} =  \sqrt{
2 \log\bigl({N \choose n} {M \choose m}\bigr)}$.

%th4.1 #&#
\begin{theorem} \label{TUV}
Assume \nref{assump1}. We suppose that alternatives under
consideration are such that
%
%e4.1 #&#
\begin{equation}
\label{altmax} %\max_{i,j} \frac{s_{ij}}{\sigma\sqrt{NM}} = o(1).
\tau_G:=\frac{\max_{(i,j) \in C}s_{ij}}{\sum_{(i,j) \in
C}s_{ij}}=\mathrm{o}(1).
\end{equation}
If the quantity $a$ is such that one of the following conditions
hold
%
%e4.2 #&#
\begin{equation}
\label{42} \frac{a}{\sigma} \sqrt{nmpq} \to\infty\quad\mbox{or}\quad
\liminf \frac{ a \sqrt{nm}}{\sqrt{2 \sigma(n \log(p^{-1})+m
\log(q^{-1}))}}>1
\end{equation}
then $\hat\psi^*$, with $H\to\infty$ such that $H^2< \mathrm{o}(a^2
nmpq/\sigma^2+\tau_G^{-1}) $ when $a \sqrt{nmpq}/\sigma\to\infty$,
% and with $T_{nm} = \sqrt{ (2+\delta) \log({N \choose n} {M \choose
%m})}$
%for some $\delta>0 $ small enough,
is such that $\gamma_{nm,a} (\psi^*)\to0$.

\end{theorem}

Proof is given in Appendix~\ref{prTUV}.

Assumption \nref{altmax} translates the fact that alternatives do
not contain too ``prominent'' values $s_{ij}$. This holds when
$\max_{(i,j)\in C}s_{ij}=\mathrm{o}(\sigma\sqrt{NM})$ and
$a\sqrt{nmpq}/\sigma\to\infty$. The assumption ensures that the
estimator of the unknown variance converges fast enough in order to
detect the signal with the same rates as in the case of known variance.
Moreover, one can slightly modify the
proof and check that the condition \nref{42} can be replaced by
the condition $\max_{(i,j)\in C}s_{ij}=\mathrm{o}(\sigma\sqrt{NM})$ when
$a\sqrt{nmpq}/\sigma=\mathrm{O}(1)$ as well.

%s4.2 #&#
\subsection{Extension to general law from an exponential family}

In many applications, we do not have Gaussian observations. Instead,
we have observations
$ X_{ij} $, i.i.d. with probability density $g_{\theta_{ij}}$
from an exponential family,
for all $ i=1,\ldots,N, $ and $j=1,\ldots,M$.
We explain here how to use the previous testing procedures in order to deal
with such setups and check that results similar to the case of Gaussian
variables hold in this case. The exponential model will behave like a
Gaussian model when the number of data is large, by asymptotic
equivalence. We expect that the optimal detection boundary is the one
for the Gaussian model properly rescaled.

We assume that the laws belong to an
exponential family in the general form
%
%e4.3 #&#
\begin{equation}
\label{genform} g_\theta(x) = \mathrm{e}^{\eta(\theta) \cdot T( x) - C(\theta)} h(x),\qquad \theta \in\Theta
\end{equation}
for the dominating measure $\mu$,
where $\eta$ is supposed 2 times continuously differentiable and strictly
increasing on $\Theta$, that is, $\eta'(\theta) > 0 $.

We consider a point $\theta^0$ interior to $\Theta$ and test, based on
$X_{ij}$'s, the null
hypothesis
$
H_0\dvt  \theta_{ij} = \theta^0 \mbox{ for all } i=1,\ldots,N,   j=1,\ldots,M,
$
against the alternative
%
%e4.4 #&#
\begin{equation}
\label{altexp} H_1\dvt\quad  \exists C \in\mathcal{C}_{nm} \mbox{
such that }\theta_{ij} = \theta^0 \mbox{ if } (i,j) \notin C \mbox{ and } \theta_{ij} - \theta^0 \geq d \mbox{ if
} (i,j) \in C.
\end{equation}

In order to build the test procedure as previously, we will rescale the
observations as follows. First, put the exponential model in the
canonical form, then change variables to $Y_{ij}=
(T(X_{ij})-m^0)/\sigma^0$, with
$m^0=E_{\theta^0} (T(X))$ and $\sigma^0=\sqrt{\operatorname{Var}_{\theta^0}(T(X))}$
computed under the null hypothesis. Let us denote the common density of
$Y_{ij}$'s by
\[
f_s(y) = \mathrm{e}^{s\cdot y - A(s)} h(y),
\]
where $s=\eta(\theta) \sigma_0 $ and $A(s)=B(s/\sigma^0) - s
m_0/\sigma_0$
and $B(\eta(\theta)) = C(\theta)$.
Here, we have $A'(s^0)=0$, $A''(s^0) = 1$ and
%
%e4.5 #&#
\begin{equation}
\label{assumpTE} A\bigl(s^0 + a\bigr) -A\bigl(s^0\bigr)
\sim\frac{a^2}2 \quad\mbox{and}\quad A'\bigl(s^0 + a
\bigr)\sim a\qquad \mbox{as } a \to0.
\end{equation}
In this way, the original problem corresponds to testing, based on
$Y_{ij}$'s, the null hypothesis
$
H_0\dvt  s_{ij} = s^0 \mbox{ for all } i=1,\ldots,N,   j=1,\ldots,M,
$
against the alternative
\[
H_1\dvt \quad \exists C \in\mathcal{C}_{nm} \mbox{ such that }
s_{ij} = s^0 \mbox{ if } (i,j) \notin C \mbox{ and }
s_{ij} - s^0 \geq a \mbox{ if } (i,j) \in C.
\]

We have the following results %general upper bounds
for exponential models.
%
%th4.2 #&#
\begin{theorem} \label{TE}
Assume \nref{assump1}. We suppose that
%
%e4.6 #&#
\begin{equation}
\label{tend} \frac{\log(p^{-1})}{m} + \frac{\log(q^{-1})}{n} \to0.
\end{equation}

\emph{Upper bounds}. If $a$ is such that one of the following
conditions hold
\[
A'\bigl(s^0+a\bigr) \sqrt{nmpq} \to\infty \quad\mbox{ or}\quad
\liminf\frac{A'(s^0+a) \sqrt{nm}}{\sqrt{2 (n \log(p^{-1})+m \log
(q^{-1}))}}>1
\]
then $\psi^*$, with $H\to\infty$ such that $H \leq c A'(s^0+a) \sqrt {nmpq} $
for some $0<c<1$ and with $T_{nm} $ replaced by $T_{nm,\delta} $
for some $\delta>0$ small enough, is such that $\gamma_{nm,a}(\psi^*)
\to0$.

\emph{Lower bounds}. Assume, moreover, that conditions \nref{cond3a}
and \nref{aa} hold. If $a$ is such that the conditions \nref
{cond12oppa} and \nref{cond12oppb} are satisfied,
then
$\g_{nm,a}\to1 $ and
$\b_{nm,a,\a}\to1-\a$ for any $\a\in(0,1)$.
\end{theorem}

Proof of the upper bounds is given in Appendix~\ref{sectionTE}.

The proof of the lower bounds uses the relation \nref{assumpTE} and
follows exactly the same lines as the proof
of Theorem~\ref{TL} in Section~\ref{PL} except that we have to consider
$T_{kl}^2 \sim(2+\delta)(k\log(p^{-1})+l \log(q^{-1}))$ for some small
$\delta>0$ instead of thresholds in \nref{Tkl}.

Under the assumption \nref{tend}, the detection boundary $a^*\to0$. Therefore,
\[
A'\bigl(s^0+a^*\bigr) \sim a^* \sim\bigl(\eta(\theta)-
\eta\bigl(\theta^0\bigr)\bigr) \sigma^0 \sim
\eta'\bigl(\theta^0\bigr) \sigma^0 d^*
\]
as $d^*\to0$. It is well known that the Fisher information at $\theta^0$ in model \nref{genform}
is $I(\theta_0) = (\sigma^0 \eta'(\theta^0))^2$. In this way, we deduce the sharp asymptotic detection boundary for
alternative~\nref{altexp} from Theorem~\ref{TE}:
$
d^* = a^* /\sqrt{I(\theta^0)}.
$

Examples of such calculations for most popular probability
distributions in the exponential family are given in Table~\ref{table2}.

%t1 #&#
\begin{table}
\caption{Examples of calculations for testing in general exponential
families}\label{table2}
\begin{tabular*}{\textwidth}{@{\extracolsep{\fill}}lllll@{}}
\hline
Probability law & $\eta$ & $m^0$ & $\sigma^0$
%& $s(\theta) = \sigma^0 \eta(\theta)$
& $\sqrt{I(\theta^0)} = \sigma^0\eta'(\theta^0)$ \\
\hline
$\operatorname{Poisson}(\theta)$, $\theta>0$& $\log(\theta)$ &
$\theta^0$ &$\sqrt {\theta^0}$
%& $\sqrt{\theta^0} \log(\theta) $
& $ (\theta^0)^{-1/2}$\\
$\operatorname{Ber}(\theta)$, $0<\theta<1$& $\log(\frac{\theta
}{1-\theta})$ & $\theta^0$& $\sqrt{\theta^0 (1-\theta^0)}$ & $(\theta^0 (1-\theta^0))^{-1/2}$
\\
$\operatorname{Exp}(\theta)$, $\theta>0$& $-\theta^{-1}$& $\theta^0 $& $\theta^0$&
$(\theta^0)^{-1}$\\
$N(0, \theta^2)$, $\theta>0$& $-\frac1{2 \theta^2}$ & $(\theta^0)^2$
& $ 2(\theta^0)^2$&
$2(\theta^0)^{-1}$ \\
\hline
\end{tabular*}
\end{table}

%s4.3 #&#
\subsection{Extension to two-sided alternative}

Let us consider model \nref{mod} and the same null hypothesis \nref
{null}, against the two-sided alternative:
\[
H_1\dvt \quad \exists C \in\mathcal{C}_{nm} \mbox{ such that
}s_{ij}=0 \mbox{ if }(i,j) \notin C \mbox{ and } |s_{ij}|
\geq a \mbox{ if } (i,j) \in C.
\]

Let us consider the following test procedures
\[
z_{\mathrm{lin}} = \frac1{\sqrt{2NM}} \sum_{i,j}
\bigl(Y_{ij}^2-1\bigr) \quad\mbox{and}\quad \psi^z_{\mathrm{lin}}
= \1_{z_{\mathrm{lin}} > H}
\]
and
\[
z_{\mathrm{max}} = \max_{C \in\mathcal{C}_{nm}} Z_C \qquad\mbox{where }
Z_C = \frac1{\sqrt{2nm}} \sum_{(i,j) \in C}
\bigl(Y_{ij}^2-1\bigr)\quad \mbox{and}\quad \psi_{\mathrm{max}}^z=
\1_{z_{\mathrm{max}} > T_{nm,\delta}}
\]
for some $\delta>0$ small enough.

%th4.3 #&#
\begin{theorem} \label{TTS}
Assume \nref{assump1}. We suppose that \nref{tend} holds.

\emph{Upper bounds}. If $a$ is such that one of the following
conditions hold
\[
a^2 \sqrt{nmpq} \to\infty\quad \mbox{or}\quad \liminf\frac{a^2 \sqrt{nm}}{2 \sqrt{ n \log(p^{-1})+m \log
(q^{-1})}} > 1
\]
then $\psi^z = \max\{\psi^z_{\mathrm{lin}}, \psi^z_{\mathrm{max}}\}$ with $H\to
\infty$
such that $H\leq c a^2/2 \sqrt{nmpq} $
for some $0<c<1$,
%and with $T_{nm} = \sqrt{ (2+\delta) \log({N \choose n} {M \choose
%m})}$ for some $\delta>0$ small enough,
is such that $\gamma_{nm,a}(\psi^z) \to0$.

\emph{Lower bounds}. Assume, moreover, that conditions \nref{cond3a}
and \nref{aa} hold. If $a$ is such that the following two conditions
are satisfied:
\[
a^2 \sqrt{nmpq} \to0 %\]
\quad\mbox{and}\quad
\limsup\frac{a^2 \sqrt{nm}}{2\sqrt{ n\log(p^{-1})+m\log(q^{-1})}}<1,
\]
then %the distinguishability is impossible, i.e.,
$\g_{nm,a}\to1 $ and
$\b_{nm,a,\a}\to1-\a$ for any $\a\in(0,1)$.
\end{theorem}
Proof is given in Appendix~\ref{ProofTS}.

%s4.4 #&#
\subsection{Related testing problems}\label{resadd}

Let us consider again the model
\nref{mod} and the null hypothesis \nref{null}. We shall see how our
alternative which locates signal in submatrices of the large matrix
compares to other alternatives. We consider first the alternatives
where the signal is located anywhere (no structure: larger alternative)
and then where the signal is located in block-submatrices (smaller alternative).

%s4.4.1 #&#
\subsubsection{Subsets without structure}
%The first problem corresponds to the alternative of the following type.
Let ${{\mathcal D}}_k$ consists of all subsets
$D\subset\{1,\ldots,N\}\times\{1,\ldots,M\}$ of cardinality
$\#(D)=k$ and let $k=nm$. Let us consider the alternative
%
%e4.7 #&#
\begin{equation}
\label{alt2} H_1\dvt\quad  \exists D\in{{\mathcal D}}_{nm}
\mbox{ such that } s_{ij}=0 \mbox{ if } (i,j)\notin D \mbox{ and }
s_{ij}\ge a \mbox{ if } (i,j)\in D
\end{equation}
(we do not suppose that the set $D$ is of product structure).
Clearly, we can consider the matrix $\{Y_{ij}\}$ as a vector of
dimension $P=NM$, and the problem is well studied as $P\to\infty$,
see Ingster~\cite{I97}, Ingster and Suslina~\cite{IS02b},
Donoho and Jin~\cite{DJ04}.

The results are as
follows. Let $k=P^{1-\b}, \b\in(0,1)$. First, let $\b\le1/2$
which corresponds to $k^2=\mathrm{O}(P)$, that is, $(nm)^2=\mathrm{O}(NM)$. Then the
detection boundary is determined by the first condition in
\nref{cond4}. It means that distinguishability is impossible when
$a^2nmpq\to0$. On the other hand, if $a^2nmpq\to\infty$, then
distinguishability
is provided by the tests of the type $\psi_H^{\mathrm{lin}}$.

Let $\b\in(1/2,1)$. Then the detection boundary is determined by
the relation
\[
a^* \sim\varphi(\b)\sqrt{\log(P)}= \varphi(\b)\sqrt{\log(NM)},
\]
where
\[
\varphi(\beta)=\cases{\sqrt{2\b-1},&\quad $1/2<\b\le3/4,$\vspace *{2pt}
\cr
\sqrt{2}(1-\sqrt{1-\b}),&\quad $3/4<\b<1,$ }\qquad \b=1-\frac{\log(nm)}{\log(NM)}.
\]
This means that, if $\limsup a /(\varphi(\b)\sqrt{\log(NM)})<1$, then
distinguishability is impossible, and if $\liminf
a/ (\varphi(\b)\sqrt{\log(NM)})>1$, then distinguishability is
provided by
the ``high criticism'' tests $\psi^{\mathrm{HC}}=\1_{\{{L_{\mathrm{HC}}>H}\}}$ based on
statistics
\[
L(t)=\frac{\sum_{i,j}(\1_{\{Y_{ij}>t\}}-\Phi(-t))}{\sqrt{NM\Phi
(t)\Phi
(-t)}},\qquad L_{\mathrm{HC}}=\max_{t_0<t}L(t),\qquad
t_0>0,
\]
with $H=\sqrt{c\log\log(NM)}, c>2.$

%s4.4.2 #&#
\subsubsection{Block-structured submatrices}
% The second problem is stated as follows.
Let ${\mathcal{E}}_{nm}$ consist of all rectangles of size $n\times
m$, that is, of the sets $E_{kl}=\{k+1,\ldots, k+n\}\times\{l+1,\ldots,
l+m\}, 0\le k\le N-n,  0\le l\le M-m$, and the alternative is of
the form
%
%e4.8 #&#
\begin{equation}
\label{alt3} H_1\dvt\quad  \exists E\in{{\mathcal E}}_{nm}
\mbox{ such that } s_{ij}=0 \mbox{ if } (i,j)\notin E \mbox{ and }
s_{ij}\ge a \mbox{ if } (i,j)\in E.
\end{equation}
Similar problems were studied recently in Arias-Castro \textit{et
al.}~\cite{ACDH} and \cite{ACCD11}
for other related geometrically-shaped clusters. Note that Arias-Castro
\textit{et al.} \cite{ACDH}
also deals with detection of rectangular shapes in a square matrix.

The detection boundary for \nref{alt3} is determined by
\[
a^*\sim\sqrt{\frac{2(\log(p^{-1})+\log(q^{-1}))}{nm}}.
\]

Let us consider the test $\psi_Z$ based on the scan statistic over a
particular set of possible rectangles, which is a suitable ``grid'' on
${\mathcal E}_{nm}$
constructed as follows.

Take $\eta_{nm}=\eta>0$. Put $n_k=(k-1)n\eta, k=1,\ldots,K, m_l=(l-1)m\eta, l=1,\ldots,L,$ where $K, L$ are such that
$N-n(1+\eta)\le n_K\le N-n, M-m(1+\eta)\le m_L\le M-m$, which yield
$K\sim N/(\eta n), L\sim M/(\eta m)$. Put
\[
Z_{kl}=\frac{1}{ \sqrt{nm}}\sum_{(i,j)\in E_{n_k m_l}}Y_{ij},\qquad
Z=\max_{1\le k\le K, 1\le l\le L}Z_{kl},\qquad  \psi_Z=\1_{Z>\sqrt{2\log(KL)}}.
\]

In this construction, we
scan over a number $K \times L$ of rectangles
which is much smaller than the cardinality of ${\mathcal E}_{nm}$ (for
technical reasons)
and which is also much larger than the set of non-overlapping rectangles
(this set would not be large enough).

%th4.4 #&#
\begin{theorem}\label{T2} $ $ Assume \nref{assump1}. Then

\emph{Upper bounds}. Let
\[
\lim\inf\frac{a^2nm}{2(\log(p^{-1})+\log(q^{-1}))}>1,
\]
and $\eta=\eta_{nm}$ is taken in such way that $\eta\to0$, $n \eta
\to
\infty$,
$ m \eta\to\infty$, $|\log(\eta)|=\mathrm{o}(|\log(pq)|)$.
Then $\g_{nm,a}(\psi_Z)\to0$ for the test
procedure $\psi_Z$ previously described.

\emph{Lower bounds}. Let
\[
\lim\sup\frac{a^2nm}{2(\log(p^{-1})+\log(q^{-1}))}<1.
\]
Then $\g_{nm,a}\to1,  \b_{nm,a,\a}\to1-\a$ for any $\a\in(0,1)$.

\end{theorem}
Proof is given in Appendix~\ref{PT2}.

Note that, the separation rates, that is, the asymptotics of $a$ that
provide distinguishability for the alternative \nref{alt}, are intermediate
between the fast separation rates for the alternative \nref{alt3} and
the slow
rates for the alternative without structure \nref{alt2}.

Let us consider the particular case of squared matrices ($N=M$) and
squared submatrices
($n=m$) such that $n=N^{1-\beta}$ for some $\beta\in(0,1)$.
The sharp asymptotic rates of the detection boundaries can be compared
in Table~\ref{table}.

%t2 #&#
\begin{table}
\caption{Table of sharp asymptotic rates of the detection boundary
$a^*$ for squared matrices and $n=N^{1-\beta}$}\label{table}
\begin{tabular*}{\textwidth}{@{\extracolsep{\fill}}llll@{}}
\hline
Rates & No structure \nref{alt2} & Submatrix \nref{alt} &
\multicolumn{1}{c@{}}{Block
structure \nref{alt3}}\\
\hline
$\beta\in(0,\frac13]$ & \multirow{2}{40pt}[-8pt]{$N^{-(1-2\beta)}$}&$N^{-(1-2\beta)}$
&\\
&&\multicolumn{1}{c}{\hrulefill}&\\
$\beta\in(\frac13, \frac12]$ &&\multirow{2}{90pt}[-8pt]{$ N^{-(1-\beta)/2} \sqrt{4 \beta
\log(N)}$}& $N^{-(1-\beta)} \sqrt{4 \beta\log
(N)}$\\
&\multicolumn{1}{c}{\hrulefill}&&\\
$\beta\in(\frac12, 1)$ & $\varphi(\beta) \sqrt{2 \log(N)} $ &&\\
\hline
\end{tabular*}
\end{table}

%%%%%%%%%%%%%%%%%%%%%%%%%%%%%%%%%%%%%%%%%%%%%%%%%%%%
%s5 #&#
\section{\texorpdfstring{Proof of Theorem \protect\ref{TL}}{Proof of Theorem 2.2}}\label{PL}
%%%%%%%%%%%%%%%%%%%%%%%%%%%%%%%%%%%%%%%%%%%%%%%%%%%%

In the first part, we give the proof of the theorem and the other parts
of this
section are dedicated to proofs of intermediate results. More lemmas are
in the \hyperref[app]{Appendix}.

We prove the lower bounds by first reducing the minimax testing error
to a Bayesian
testing risk with uniform prior over the set of parameters.
Typically, one studies the likelihood ratio under the prior with
respect to the law
$P_0$ under the null hypothesis and proves that it tends to 1
in quadratic mean (under $P_0$). Nevertheless, this does not work as the
covariance of the likelihood ratio is too large. Therefore, we truncate the
likelihood ratio in a convenient way.

%s5.1 #&#
\subsection{Prior and truncated likelihood ratio}

%In order to simplify the notation, we consider the case $\sigma=1$
%only.
Let $S_C=\{s_{ij}\}$ be the matrix such that $s_{ij}=0, (i,j)\notin C,  s_{ij}=a, (i,j)\in C$. Let us consider the prior
on the set of matrices:
\[
\pi=G_{nm}^{-1}\sum_{C\in{\mathcal C}_{nm}}
\delta_{S_C},\qquad G_{nm}=\# ({\mathcal C}_{nm}),
\]
and let $P_\pi$ be the mixture of likelihoods
$ P_\pi=G_{nm}^{-1}\sum_{C\in{\mathcal C}_{nm}}P_{S_C}. $ Let us
consider the
likelihood ratio
\[
L_\pi(Y)=\frac{\mathrm{d}P_\pi}{\mathrm{d}P_0}(Y)=G_{nm}^{-1}
\sum_{C\in{\mathcal
C}_{nm}}\frac{\mathrm{d}P_{S_C}}{\mathrm{d}P_0}(Y)=G_{nm}^{-1}
\sum_{C\in{\mathcal
C}_{nm}}\exp\bigl(-b^2/2+bY_C
\bigr),
\]
here and below we set $b^2\eq a^2nm$, and, for submatrix $C$ of the
size $n\times m$, the statistics $Y_C$ are defined by \nref{Y}.
Since $\pi({\mathcal S}_{nm})=1$, in order to obtain
indistinguishability: $\g_{nm,a}\to1, \b_{nm,a,\a}\to1-\a, \forall  \a\in(0,1)$, it suffices to show
%
%e5.1 #&#
\begin{equation}
\label{LR1} L_\pi(Y)\to1\qquad  \mbox{in } P_0\mbox{-probability}.
\end{equation}
Indeed,
\begin{eqnarray*}
\g_{nm,a} & = & \inf_{\psi\in[0,1]} \sup_{S \in{\mathcal S}_{nm,a}} \bigl( \a(
\psi) + \beta(\psi,S) \bigr)
\\
& \geq&\inf_{\psi\in[0,1]} \frac1{G_{nm}} \sum
_{S \in{\mathcal S}_{nm,a}} \biggl( E_0\bigl(\psi(Y)\bigr) +
E_0 \biggl[ \bigl(1-\psi(Y)\bigr) \frac
{\mathrm{d}P_{S_C}}{\mathrm{d}P_0}(Y) \biggr]
\biggr)
\\
& \geq&\inf_{\psi\in[0,1]} \bigl( E_0\bigl(\psi(Y)\bigr) +
E_0 \bigl[ \bigl(1-\psi(Y)\bigr) L_\pi(Y) \bigr] \bigr)
\\
& \geq& E_0\bigl(\psi^*(Y)\bigr) + E_0 \bigl[ \bigl(1-
\psi^*(Y)\bigr) L_\pi(Y) \bigr],
\end{eqnarray*}
where $\psi^*(Y) = \1_{L_\pi(Y) >1}$ is the likelihood ratio test. Therefore,
\nref{LR1} implies by Fatou's lemma that
\[
\liminf\g_{nm,a} \geq E_0 \bigl[ \liminf \bigl( \psi^*(Y) +
\bigl(1-\psi^*(Y)\bigr) L_\pi(Y) \bigr) \bigr],
\]
that is, $\g_{nm,a} \to1$. It is easy to deduce that
$\beta_{nm,a,\a} \to1-\a.$

Let us replace the statistics $L_\pi(Y)$ by their truncated version
\[
\tilde L_\pi(Y)=G_{nm}^{-1}\sum
_{C\in{\mathcal C}_{nm}} \frac{
\mathrm{d}P_{S_C}}{\mathrm{d}P_0}(Y)\1_{\Gamma_C},
\]
where the events $\Gamma_C$ are determined as follows. Set
\[
T_{kl}=\sqrt{2\bigl(\log(G_{kl})+\log(nm)\bigr)}\to\infty.
\]
Take small $\delta_1>0$ (which will be specified later) and set
$k_0=\delta_1 n, l_0=\delta_1 m$. Let $\CC_{kl,C}=\{V\in\CC_{kl}\dvt  V \subset C\}$ be the submatrices of $C \in\CC_{nm}$ which are in
$\CC_{kl}$.
Then we set
%
%e5.2 #&#
\begin{equation}
\label{ZC} \Gamma_C=\bigcap_{k_0\le k\le n, l_0\le l\le m}
\bigcap_{V \in\CC_{kl,C}}\{Y_V \le T_{kl}
\}.
\end{equation}
By \nref{L1}, under conditions on $k,l$ in \nref{ZC} (and similarly to
the equivalent of $T^2_{nm}$) we have
%
%e5.3 #&#
\begin{equation}
\label{Tkl} T_{kl}^2\sim2\bigl(k\log\bigl(p^{-1}
\bigr)+l\log\bigl(q^{-1}\bigr)\bigr).
\end{equation}

Indeed, when looking at second-order moments of the likelihood ratio
$L_\pi(Y)$ a large contribution comes from overlapping submatrices
$C_1$ and $C_2$
inducing correlated random variables $Y_{C_1}$ and $Y_{C_2}$.
Our idea is to truncate $Y_V$, for submatrices $V$ of size close to $(k,l)$,
at its expected maximal value in order
to reduce the contribution of these correlations.
%hence the likelihood $\tilde L_\pi(Y)$.

%pr5.1 #&#
\begin{proposition}\label{PZ} Set $\Gamma_{nm}=\bigcap_{C\in
\CC_{nm}}\Gamma_C$. Then, under the assumptions of Theorem~\ref{TL}, $
P_0(\Gamma_{nm})\to1. $
\end{proposition}
Proof is given in Appendix~\ref{PPZ}.

Proposition \ref{PZ} yields
\[
P_0 \bigl(L_\pi(Y)=\tilde
L_\pi(Y) \bigr)\to1,
\]
and in place of \nref{LR1} it suffices to check that
%
%e5.4 #&#
\begin{equation}
\label{LR4} \tilde L_\pi(Y)\to1 \qquad\mbox{in } P_0\mbox{-probability}.
\end{equation}

In order to get \nref{LR4} it suffices to verify two relations:
%
%pr5.2 #&#
\begin{proposition}\label{PE} Under the assumptions of Theorem~\ref
{TL}, we have
\[
E_0(\tilde L_\pi)\to1.%,  E_0(\tilde L_\pi^2)\le1+o(1).
\]
\end{proposition}
Proof is given in Appendix~\ref{PPE}.

%pr5.3 #&#
\begin{proposition}\label{PU} Under the assumptions of Theorem~\ref
{TL}, we have
\[
E_0\bigl(\tilde L_\pi^2\bigr)\le1+\mathrm{o}(1).
\]
\end{proposition}

Propositions~\ref{PE} and~\ref{PU} imply that
\[
E_0(\tilde L_\pi-1)^2= \bigl(E_0
\bigl(\tilde L_\pi^2\bigr)-1 \bigr)-2
\bigl(E_0(\tilde L_\pi)-1 \bigr)\le \mathrm{o}(1),
\]
which ends the proof of the theorem.

The remaining part of this section is devoted to obtaining the
Proposition \ref{PU}.

%s5.2 #&#
\subsection{\texorpdfstring{Proof of Proposition \protect\ref{PU}}{Proof of Proposition 5.3}}\label{SOM}

We deal with the second order moment of the truncated likelihood ratio.
We have
%
%e5.5 #&#
\begin{equation}
\label{LR7} E_0\bigl(\tilde L_\pi^2
\bigr)=G_{nm}^{-2}\sum_{C_1\in{\mathcal C}_{nm},C_2\in
{\mathcal
C}_{nm}}E_0
\bigl(\exp\bigl(-b^2+b(Y_{C_1}+Y_{C_2})\bigr)
\1_{\{\Gamma
_{C_1}\cap
\Gamma_{C_2}\}} \bigr).
\end{equation}
We note that the expected value in the previous sum does not depend on
$C_1$ and $C_2$ but merely on the size of their common submatrix.
Let $C_1=A_1\times B_1, C_2=A_2\times B_2$ and set
\[
k=\#(A_1\cap A_2),\qquad l=\#(B_1\cap
B_2),\qquad V=(A_1\cap A_2)\times(B_1
\cap B_2).
\]
Under this notation, we put
\[
g(k,l)=E_0 \bigl(\exp\bigl(-b^2+b(Y_{C_1}+Y_{C_2})
\bigr)\1_{\{\Gamma_{C_1}\cap
\Gamma_{C_2}\}} \bigr)
\]
and see that we can rewrite \nref{LR7} as follows:
\[
E_0\bigl(\tilde L_\pi^2\bigr) = \sum
_{k=0}^n \sum
_{l=0}^m \frac{\#((C_1,C_2)\in\mathcal{C}^2_{nm}\dvt \mbox{size}(V)=(k,l))}{G_{nm}^2} g(k,l) .
\]
%
%In view of symmetry we can fix $C_1$, say $C_1=\{1,\ldots n\}\times
%E_0(\tilde L_\pi^2)=G_{nm}^{-1}\sum_{C_2\in{\mathcal C}_{nm}}g(k,l).
%%\le
%%G^{-1}\sum_{C_2\in{\mathcal C}_{nm}}g(k,l).

\textit{Notation}.
Set $z_{kl}^2=a^2kl, \rho_{kl}=kl/nm$ and recall that $b^2=a^2 n m$.
Note that it means also that $b^2 = z^2_{mn}$ and that $z^2_{kl} = b^2
\rho_{kl}$.

%le5.1 #&#
\begin{lemma}\label{LGen}
The following inequalities hold true.
\begin{longlist}[(1)]
\item[(1)] We have
%
%e5.6 #&#
\begin{equation}
\label{in1} g(k,l)\le E_0 \bigl(\exp\bigl(-b^2+b(Y_{C_1}+Y_{C_2})
\bigr) \bigr)=\exp \bigl(z^2_{kl}\bigr)\eq
g_1(k,l).
\end{equation}

\item[(2)] Let $b\ge T_{nm}/(1+\rho_{kl})$. Then
%
%e5.7 #&#
\begin{eqnarray}\label{in2}
g(k,l)&\le& E_0 \bigl(\exp\bigl(-b^2+b(Y_{C_1}+Y_{C_2})
\bigr) \1_{\{Y_{C_1}\le
T_{nm},  Y_{C_2}\le T_{nm}\}} \bigr)
\nonumber
\\[-8pt]
\\[-8pt]
\nonumber
&\le& \exp \biggl(-(T_{nm}-b)^2+
\frac{\rho_{kl}T_{nm}^2}{1+\rho_{kl}} \biggr) \eq g_2(k,l).
\end{eqnarray}

\item[(3)] Let $k\ge\delta_1 n, l\ge\delta_1 m$, and $T_{kl}\le
2z_{kl}$. Then
%
%e5.8 #&#
\begin{eqnarray}\label{in3}
g(k,l)&\le& E_0 \bigl(\exp\bigl(-b^2+b(Y_{C_1}+Y_{C_2})
\bigr) \1_{\{Y_V \le
T_{kl}\}} \bigr)
\nonumber
\\[-8pt]
\\[-8pt]
\nonumber
&=&\exp\bigl(T^2_{kl}/2-(T_{kl}-z_{kl})^2
\bigr)\eq g_3(k,l).
\end{eqnarray}
\end{longlist}
\end{lemma}

Proof of Lemma \ref{LGen} is given in Appendix~\ref{LGenP}.

%s5.2.1 #&#
\subsubsection{From hypergeometric to binomial distributions}

Observe that the right-hand side of \nref{LR7} is the expectation of
$g(X_1,X_2)$ over $X_1,X_2$ which are independent and having
hypergeometric distributions ${{\mathcal{HG}}}_1={{\mathcal{HG}}}(N,n,n),
{{\mathcal{HG}}}_2={{\mathcal{HG}}}(M,m,m)$, respectively, that is,
%
%e5.9 #&#
\begin{equation}
\label{HG1} E_0\bigl(\tilde L_\pi^2\bigr)=
\sum_{k=0}^n \sum
_{l=0}^m \frac{{N\choose n} {n \choose k}{{N-n} \choose{n-k}}\cdot
{M\choose m} {m \choose l}{{M-m} \choose{m-l}}} {
{N \choose n}^2 \cdot{M \choose m}^2} g(k,l)
=E_{{{\mathcal{HG}}}_1\times{{\mathcal{HG}}}_2}g(X_1,X_2).
\end{equation}

Let us compare random variables $X$ having hypergeometric
distributions ${{\mathcal{HG}}}={{\mathcal{HG}}}(N,n,n)$ and binomial
distribution $\mathrm{Bin}=\mathrm{Bin}(n, \tilde p), \tilde p=n/(N-n)$.

%le5.2 #&#
\begin{lemma}\label{L4a}
Under binomial distribution, $X$ is stochastically larger, than
under hypergeometric distributions, that is, for any $x\in\R$,
\[
P_{{{\mathcal{HG}}}}(X\ge x)\le P_{\mathrm{Bin}}(X\ge x).
\]
This yields, for any non-decreasing function $g$,
\[
E_{{{\mathcal{HG}}}}\bigl(g(X)\bigr)\le E_{{{\mathrm{Bin}}}}\bigl(g(X)\bigr).
\]
\end{lemma}
\begin{pf} The first claim corresponds to Lemma 3 in Arias-Castro
\textit{et al.} \cite{arias}.
The second claim follows from the Abel's transform of the series for
the expectation.
\end{pf}

Let $P_{n,p}(k)=P_{\mathrm{Bin}}(X=k)$, for some integer $k$, where $X$ has
binomial $\mathrm{Bin}(n,p)$ distribution, and similarly
$P_{N,n,n}(k)=P_{{{\mathcal{HG}}}}(X=k)$ for hypergeometric distributions
${{\mathcal{HG}}}(N,n,n)$ of $X$. %We use the following
%
%le5.3 #&#
\begin{lemma}\label{L5} Let $n\to\infty, p\to0,  p>0, k\ge
n/r(p)$ where $r(p)\ge1$ for $p>0$ small enough, and
$\log(r(p))=\mathrm{o}(\log(p^{-1}))$. Then
\begin{eqnarray*}
\log\bigl(P_{n,p}(k)\bigr)&\le& k\log(p) \bigl(1+\mathrm{o}(1)\bigr),
\\
\log\bigl(P_{N,n,n}(k)\bigr)&\le&k\log(p) \bigl(1+\mathrm{o}(1)\bigr).
\end{eqnarray*}
\end{lemma}
Proof is given in Appendix~\ref{PL5}.

%s5.2.2 #&#
\subsubsection{Evaluation of the expectation}

Take any small $\delta>0$. The detection boundary $a$ satisfies assumption
\nref{cond12oppb}, where the worst-case is when the limit is close to 1.
It suffices, therefore, to consider the case
%
%e5.10 #&#
\begin{equation}
\label{b} b^2=a^2nm\sim(2-\delta) \bigl(n\log
\bigl(p^{-1}\bigr)+m\log\bigl(q^{-1}\bigr)\bigr).
\end{equation}
This implies
%
%e5.11 #&#
\begin{equation}
\label{a} a^2\asymp\frac{\log(p^{-1})}{m}+\frac{\log(q^{-1})}{n}.
\end{equation}

In order to evaluate the right-hand side of \nref{HG1}, let us
firstly divide the expectation into two parts $E_{{{\mathcal{HG}}}_1\times{{\mathcal{HG}}}_2}
(g(X_1,X_2))=E_1+E_2$, where
\[
E_1=E_{{{\mathcal{HG}}}_1\times{{\mathcal{HG}}}_2} \bigl(g(X_1,X_2)
\1_{X_1a^2<1}\bigr),\qquad E_2=E_{{{\mathcal{HG}}}_1\times
{{\mathcal{HG}}}_2} \bigl(g(X_1,X_2)
\1_{X_1a^2\ge1}\bigr).
\]

Recall that we denote $\tilde p=n/(N-n)$, $\tilde q=m/(M-m)$ and the
binomial distributions $\mathrm{Bin}_1 = \mathrm{Bin}(n,\tilde p)$ and $\mathrm{Bin}_2 =
\mathrm{Bin}(m,\tilde q)$.

We would like to show that $E_1\le1+\mathrm{o}(1)$ and $E_2=\mathrm{o}(1)$, so we keep
in mind from now on that $N,  M$ are sufficiently large.

\textit{Evaluation of $E_1$}.

It follows from \nref{a} and \nref{aa} that
$X_1=\mathrm{O}(n/\log(q^{-1}))$ under $a^2X_1< 1$. By \nref{in1} we have
\begin{eqnarray*}
E_1&\le& E_{{{\mathcal{HG}}}_1\times{{\mathcal{HG}}}_2}\bigl(\exp\bigl(a^2X_1X_2
\bigr)\1_{X_1a^2< 1}\bigr)\\&=&E_{{{\mathcal{HG}}}_1} \bigl(E_{{{\mathcal{HG}}}_2} \bigl(
\exp\bigl(a^2X_1X_2\bigr) \bigr)
\1_{X_1a^2< 1} \bigr).
\end{eqnarray*}
In view of Lemma \ref{L4a}, the expected value of a non-decreasing
function of $X_2$
having hypergeometric distribution ${\mathcal{HG}}_2$ is less than the
same expected value under the binomial $\mathrm{Bin}_2=\mathrm{Bin}(m,\tilde q), \tilde q=m/(M-m)$,
\begin{eqnarray*}
&&E_{{{\mathcal{HG}}}_2} \bigl(\exp\bigl(a^2X_1X_2
\bigr) \bigr)\\
&&\quad\le E_{\mathrm{Bin}_2} \bigl(\exp\bigl(a^2X_1X_2
\bigr) \bigr)=\bigl(1+\tilde q\bigl(\mathrm{e}^{a^2X_1}-1\bigr)\bigr)^m
\\
&&\quad\le \exp\bigl(m\tilde q\bigl(\mathrm{e}^{a^2X_1}-1\bigr)\bigr).
\end{eqnarray*}
Observe that, for some $B>0$ under constraint $X_1a^2\le1$,
\[
%E_l\l(\exp(a^2kl))\r)=(1+q(\mathrm{e}^{a^2k}-1))^m\le\exp(mq(\mathrm{e}^{a^2k}-1))
\exp\bigl(m\tilde q\bigl(\mathrm{e}^{a^2X_1}-1\bigr)\bigr)\le\exp
\bigl(Bmqa^2X_1\bigr).
\]
Taking the expectation over $X_1$, we get similarly
\begin{eqnarray*}
E_1&\le& E_{{{\mathcal{HG}}}_1}\bigl(\exp\bigl(Bmqa^2X_1
\bigr)\bigr)
\\[-1pt]
&\le& E_{\mathrm{Bin}_1}\bigl(\exp\bigl(B m q a^2X_1
\bigr)\bigr)= \bigl(1+\tilde p\bigl(\mathrm{e}^{Bm\tilde q
a^2}-1\bigr)\bigr)^n
\\[-1pt]
&\le& \exp\bigl(n\tilde p\bigl(\mathrm{e}^{Bmq a^2}-1\bigr)\bigr).
\end{eqnarray*}
By \nref{a} and \nref{aa}, we have $a^2m\asymp\log(p^{-1})$. By
condition \nref{cond3a}, we have $q\log(p^{-1})=\mathrm{o}(1)$, which yields
$mqa^2=\mathrm{o}(1)$. Thus, $n\tilde p(\mathrm{e}^{Bmqa^2}-1)\sim Bnpmqa^2=\mathrm{o}(1)$ by
the condition \nref{cond12oppa}, and we get
\[
E_1\le\exp\bigl(\mathrm{o}(1)\bigr)=1+\mathrm{o}(1).
\]

\textit{Evaluation of $E_{2}$}.

In order to evaluate $E_2$, we use the assumption \nref{aa} and write
\[
a^2 \asymp\log\bigl(p^{-1}\bigr)/m \asymp\log
\bigl(q^{-1}\bigr)/n
\]
instead of \nref{a}. Therefore, we can take $\delta_1>0$ small enough
such that
$\delta_1 a^2m\le\log(p^{-1})/2$ and $\delta_1 a^2n \le
\log(q^{-1})/2$ for $N,  M,   n $ and $m$ large enough.

Divide $E_2$ into two parts $E_2=E_{21}+E_{22}$, where
\begin{eqnarray*}
E_{21}&=&E_{{{\mathcal{HG}}}_1\times{{\mathcal{HG}}}_2} \bigl(g(X_1,X_2)
\1_{X_1a^2\ge1, X_2/m<\delta_1}\bigr),
\\[-1pt]
E_{22}&=&E_{{{\mathcal{HG}}}_1\times{{\mathcal{HG}}}_2} \bigl(g(X_1,X_2)
\1_{X_1a^2\ge1, X_2/m\ge\delta_1}\bigr).
\end{eqnarray*}

\textit{Evaluation of $E_{21}$}.

By applying Lemma \ref{L5} with $\log r(p)\eq\log(\log(q^{-1}))$,
which is $\mathrm{o}(\log(p^{-1}))$
by \nref{cond3a}, and since $P_{M,m,m}(l)\le1$, we get
\begin{eqnarray*}
E_{21}&\le& \sum_{n\ge k>a^{-2},   0\le l\le\delta_1 m} \exp
\bigl(a^2kl\bigr)P_{N,n,n}(k)P_{M,m,m}(l)
\\[-1pt]
&\le& \sum_{n\ge k>a^{-2},   0\le l\le\delta_1
m}\exp \bigl(k \bigl(a^2l-
\log\bigl(p^{-1}\bigr) \bigl(1+\mathrm{o}(1)\bigr) \bigr) \bigr).
\end{eqnarray*}
Observe that under the constraints in the sum,
\[
a^2l\le\delta_1 a^2m\le\log
\bigl(p^{-1}\bigr) \bigl(1/2+\mathrm{o}(1)\bigr),
\]
which yields in the previous exponential
\begin{eqnarray*}
k \bigl(a^2l-\log\bigl(p^{-1}\bigr) \bigl(1+\mathrm{o}(1)\bigr)
\bigr)&\le& -k\log\bigl(p^{-1}\bigr) \bigl(1/2+\mathrm{o}(1)\bigr)
\\[-1pt]
&\le& - a^{-2}\log\bigl(p^{-1}\bigr) \bigl(1/2+\mathrm{o}(1)\bigr)
\asymp m.
\end{eqnarray*}
Therefore, we have $ E_{21}\le nm\exp(-Bm)=\mathrm{o}(1) $ for some $B>0$ by
condition \nref{cond3a}.\vadjust{\goodbreak}

\textit{Evaluation of $E_{22}$}.

In order to evaluate the item $E_{22}$, we divide it in two parts as
well: $E_{22}=I_1+I_2$,
\begin{eqnarray*}
I_1&=&E_{{{\mathcal{HG}}}_1\times{{\mathcal{HG}}}_2} \bigl(g(X_1,X_2)
\1_{X_1a^2\ge
1, X_1/n< \delta_1,  X_2/m\ge\delta_1}\bigr),
\\
I_2&=&E_{{{\mathcal{HG}}}_1\times{{\mathcal{HG}}}_2} \bigl(g(X_1,X_2)
\1_{X_1/n\ge
\delta_1,  X_2/m\ge\delta_1}\bigr).
\end{eqnarray*}
The evaluation of $I_1$ is similar to the evaluation of $E_{21}$ and
we get $I_1=\mathrm{o}(1)$.

\textit{Evaluation of $I_{2}$}.

Let us divide the set $\CH= \{(k,l)\dvt  k/n \geq\delta_1,   l/m \geq
\delta_1 \}$, appearing in $I_2$, into two parts:
\begin{eqnarray*}
\CH_1&=& \bigl\{(k,l)\in\CH\dvt  k\log\bigl(p^{-1}\bigr)+l\log
\bigl(q^{-1}\bigr)\ge 2\rho_{kl}\bigl(n\log
\bigl(p^{-1}\bigr)+m\log\bigl(q^{-1}\bigr)\bigr) \bigr\},
\\
\CH_2&=& \bigl\{(k,l)\in\CH\dvt  k\log\bigl(p^{-1}\bigr)+l\log
\bigl(q^{-1}\bigr)< 2\rho_{kl}\bigl(n\log\bigl(p^{-1}
\bigr)+m\log\bigl(q^{-1}\bigr)\bigr) \bigr\}.
\end{eqnarray*}
This yields the division of $I_2$ into $I_2=I_{12}+I_{22}$. Observe
that $\rho_{kl}\ge\delta_1^2$ for $(k,l)\in\CH$.

%First let the condition \nref{cond12oppb} hold true.
Let us consider $I_{12}$. Recalling \nref{in2},
observe that we can take $\delta>0$ small enough in \nref{b} such
that $ t=T_{nm}-b(1+\rho_{kl})<0. $
%Note that if $t<0$, then the
%power in the exponent in \nref{fg} can be rewritten in the form
%$$
%b^2\rho_{kl}-\frac{(b(1+\rho_{kl})-T_{nm})^2}{1+\rho_{kl}}=
%-(T_{nm}-b)^2+\frac{T_{nm}^2\rho_{kl}}{1+\rho_{kl}},
%$$
Applying \nref{in2} and Lemma \ref{L5}
for
$P_{N,n,n}(k), P_{M,m,m}(l)$, we get
\[
I_{12}\le\sum_{(k,l)\in\CH_1}\exp
\biggl(-(T_{nm}-b)^2 +\frac{T_{nm}^2\rho_{kl}}{1+\rho_{kl}}-k\log
\bigl(p^{-1}\bigr)-l\log \bigl(q^{-1}\bigr)+\mathrm{o}
\bigl(T_{nm}^2\bigr) \biggr).
\]
Note that for $\delta>0$ small enough in \nref{b} one can take
$\delta_2=\delta_2(\delta)>0$ such that $(T_{nm}-b)^2\ge\delta_2
T_{nm}^2$ for the first item in the exponent. Denote
\[
A=A_{n,p}=n\log\bigl(p^{-1}\bigr);\qquad B=B_{m,q}= m
\log\bigl(q^{-1}\bigr)
\]
and $T_{nm}^2 = 2 \log(G_{nm}) \sim2(A+B)$ by \nref{L1}. Observe
that $2(A+B)\rho_{kl}\le A\frac kn+B\frac lm$ for $(k,l)\in
\mathcal{H}_1$.

The following terms in the power of the exponential
above can be bounded on the set $\mathcal{H}_1$ as follows
\begin{eqnarray*}
\frac{T_{nm}^2\rho_{kl}}{1+\rho_{kl}}-k\log\bigl(p^{-1}\bigr)-l\log\bigl(q^{-1}
\bigr) &= &\frac{2(A+B)\rho_{kl}}{1+\rho_{kl}}-\biggl(A\frac kn+B\frac lm\biggr)+\mathrm{o}
\bigl(T_{nm}^2\bigr)
\\
&\leq& \biggl(\frac1{1+\rho_{kl}}-1\biggr) \biggl(A\frac kn+B\frac lm
\biggr)+\mathrm{o}\bigl(T_{nm}^2\bigr)\leq \mathrm{o}\bigl(T_{nm}^2
\bigr).
\end{eqnarray*}
Therefore,
\[
I_{12}\le2nm\exp\bigl(-\bigl(\delta_2+\mathrm{o}(1)\bigr)
T_{nm}^2\bigr)=\mathrm{o}(1).
\]

Consider now the item $I_{22}$. Recalling \nref{Tkl}, \nref{b} and
$z_{kl}^2=\rho_{kl}a_{nm}^2$ observe that the constraint in $\CH_2$
corresponds to $T_{kl}^2<2z_{kl}^2(1+\mathrm{o}(1))$. This implies
$T_{kl}^2<4z_{kl}^2$ and $T_{kl}-2z_{kl}<0$ for $N, M$ large enough.

Applying \nref{in3} and Lemma \ref{L5} for $P_{N,n,n}(k), P_{M,m,m}(l)$, we similarly get
\[
I_{22}\le\sum_{(k,l)\in
\CH_2}\exp
\bigl(T_{kl}^2/2-(T_{kl}-z_{kl})^2-
\bigl(k\log\bigl(p^{-1}\bigr)+l\log \bigl(q^{-1}\bigr)\bigr)
\bigl(1+\mathrm{o}(1)\bigr) \bigr).
\]
Since $k\log(p^{-1})+l\log(q^{-1})\sim T_{kl}^2/2$, the power in
the exponent is of the form
\[
-(T_{kl}-z_{kl})^2+\mathrm{o}\bigl(T_{kl}^2
\bigr).
\]
Under \nref{b}, %(compare with \nref{cond*}) recalling \nref{ZT}
we can see that $(T_{kl}-z_{kl})^2\ge\delta_2T_{kl}^2$ for some
$\delta_2>0$ and $N, M$ large enough. In fact, recalling $A>0, B>0, k/n\in(\delta_1,1],  l/m\in(\delta_1,1]$ before we have
\begin{eqnarray*}
T_{kl}^2-z_{kl}^2&=&2\biggl(A\frac
kn+B\frac lm\biggr)-(2-\delta)\rho_{kl}(A+B)+\mathrm{o}\bigl(T_{kl}^2
\bigr)
\\
&=&\delta(A+B)\rho_{kl}+2A\frac kn\biggl(1-\frac kn\biggr)+2B\frac lm
\biggl(1-\frac lm\biggr)+\mathrm{o}\bigl(T_{kl}^2\bigr)
\\
&\ge& \delta(A+B)\rho_{kl}+\mathrm{o}\bigl(T_{kl}^2\bigr)
\sim\delta z_{kl}^2.
\end{eqnarray*}
Since $T_{kl}\asymp T_{nm}$ for $(k,l)\in\mathcal{H}$, these yield
$I_{22}=\mathrm{o}(1)$.

Proposition \ref{PU} follows.

\begin{appendix}\label{app}
%%%%%%%%%%%%%%%%%%%%%%%%
%s6 #&#
\section*{Appendix}

\setcounter{subsection}{0}
%s6.1 #&#
\subsection{\texorpdfstring{Proof of Theorem~\protect\ref{TU}}
{Proof of Theorem 2.1}}\label{PPmax}
It is easy to see that $\alpha(\psi^*) \leq\alpha(\psi_H^{\mathrm{lin}})+\alpha
(\psi^{\mathrm{max}})$ and that
$\beta_{nm,a}(\psi^{*})\leq\min\{
\beta_{nm,a}(\psi_H^{\mathrm{lin}}),\break
\beta_{nm,a}(\psi^{\mathrm{max}}) \}$.
Therefore, we study the two test procedures separately.

We have, for any real number $H$,
\[
\a\bigl(\psi^{\mathrm{lin}}_H\bigr)=\Phi(-H),\qquad \b_{nm,a}
\bigl(\psi^{\mathrm{lin}}_H\bigr)\le \Phi(H-a\sqrt{nmpq}),
\]
where $\Phi$ denotes the cumulative distribution function of a
standard Gaussian
random variable.
Indeed, observe that the statistic $t_{\mathrm{lin}}$ is standard Gaussian under
$P_0$ which yields the relation for $\a(\psi^{\mathrm{lin}}_H)$. Also
$t_{\mathrm{lin}}\sim\CN(h_{S_C}, 1)$ under $P_{S_C}$-probability,
$S_C\in\CS_{nm,a}$, where
\[
h_{S_C}=\frac{1}{\sqrt{NM}}\sum_{(i,j)\in C}s_{ij}
\ge \frac{anm}{\sqrt{NM}}=a\sqrt{mnpq}.
\]
This yields the relation $\b(\psi^{\mathrm{lin}}_H, S)\le\Phi(H-a\sqrt {mnpq}), \forall S\in\CS_{nm,a}$ and the same inequality for
$\b_{nm,a}(\psi^{\mathrm{lin}}_H)$.

As a consequence, if $H\to\infty$ then
$\a(\psi_H^{\mathrm{lin}}) = \Phi(-H) \to0$ and if \nref{cond1} holds and $H
\leq ca\sqrt{nmpq}$,
\[
\beta_{nm,a}\bigl(\psi_H^{\mathrm{lin}}\bigr) \leq\Phi
\bigl((c-1)a\sqrt{nmpq} \bigr)\to0
\]
for $c<1$. Thus $\gamma_{nm,a}(\psi^{\mathrm{lin}}_H) \to0$.\vadjust{\goodbreak}

Now, we prove that,
$
\a(\psi^{\mathrm{max}})\to0 %, \a(\psi^{\mathrm{max}})\to0,
$
and
$
\b_{nm,a}(\psi^{\mathrm{max}})\le\Phi(T_{nm}-a\sqrt{nm})
$, under \nref{assump1}.
It is not hard to
check that under \nref{assump1},
%
%e6.1 #&#
\setcounter{equation}{0}
\begin{equation}
\label{L1} \log(G_{nm})\sim\bigl(n\log\bigl(p^{-1}\bigr)+m
\log\bigl(q^{-1}\bigr)\bigr).
\end{equation}
Observe that $Y_C\sim\CN(0,1)$ under $P_0$ and, since
$G_{nm}\to\infty$ and $\Phi(-T)\asymp\exp(-T^2/2)/T$ as
$T\to\infty$, we get
\[
\a\bigl(\psi^{\mathrm{max}}\bigr)=P_0(t_{\mathrm{max}}>T_{nm})
\le\sum_{C\in\CC
_{nm}}P_0(Y_C>T_{nm})=G_{nm}
\Phi(-T_{nm})\to 0.
\]

Let $S_C\in\CS_{nm,a}$. Then $Y_C\sim\CN(g_{S_C}, 1)$ under
$P_{S_C}$-probability with
\[
g_{S_C}=(nm)^{-1/2}\sum_{(i,j)\in C}s_{ij}
\ge a\sqrt{nm}.
\]
Observe that
\begin{eqnarray*}
\b\bigl(\psi^{\mathrm{max}},S_C\bigr)&=&P_{S_C}(t_{\mathrm{max}}
\le T_{nm})\le P_{S_C}(Y_C\le T_{nm})=
\Phi(T_{nm}-g_{S_C})
\\
&\le& \Phi(T_{nm}-a\sqrt{nm}).
\end{eqnarray*}
Thus, under \nref{cond2},
\[
a\sqrt{nm} - T_{nm} = T_{nm} \biggl( \frac{a\sqrt{nm}(1+\mathrm{o}(1))}{ \sqrt{2 (n \log(p^{-1})
+ m \log(q^{-1}))}}- 1
\biggr) \to\infty
\]
and this implies that $\gamma_{nm,a}(\psi^{\mathrm{max}}) \to0$.

%%%%%%%%%%%%%%%%%%%%%%%%

%s6.2 #&#
\subsection{\texorpdfstring{Proof of Theorem \protect\ref{TLa}}
{Proof of Theorem 2.3}}\label{proofTLa}

Note that the test $\psi^{\mathrm{lin}}_H$ does not depend on $n,m$.
Therefore for distinguishability in the adaptive problem it is
sufficient to assume
\nref{cond1a} (which is a uniform version of \nref{cond1}). We have
$\gamma_{NM,\mathbf{a}}(\psi_H^{\mathrm{lin}}) \to0$.

For the test $\psi_{NM}^{\mathrm{max}},$ we obtain similarly to the nonadaptive
case, that
\begin{eqnarray*}
\a\bigl(\psi^{\mathrm{max}}_{NM}\bigr)&=&P_0(t_{NM,\mathrm{max}}>1)
\le \sum_{n,m}\sum_{C\in\CC_{nm}}P_0(Y_C>H_{nm})
=\sum_{n,m}\sum_{C\in\CC_{nm}}
\Phi(-H_{nm})
\\
&\asymp& \sum_{n,m}\sum
_{C\in\CC_{nm}}\frac{1}{NMG_{nm} \sqrt{\log(NMG_{nm})}} =\mathrm{O} \biggl(\frac1{\sqrt{\log(NM)}}
\biggr)\to 0.
\end{eqnarray*}
Moreover,
\begin{eqnarray*}
\b\bigl(\psi^{\mathrm{max}}_{NM},S_C
\bigr)&=&P_{S_C}(t_{NM,\mathrm{max}}\le1)\le P_{S_C}(Y_C
\le V_{nm})=\Phi(V_{nm}-g_{S_C})
\\
&\le& \Phi(V_{nm}-a_{nm}\sqrt{nm})
\end{eqnarray*}
and we deduce that
\[
\beta_{nm,\mathbf{a}}\bigl(\psi_{NM}^{\mathrm{max}}\bigr) \leq\Phi
\Bigl(\max_{(n,m) \in{\mathcal K}_{NM}} (V_{nm}-a_{nm}\sqrt{nm})\Bigr).
\]
By \nref{L1}, we have
\begin{eqnarray*}
&&\min_{(n,m) \in{\mathcal{K}}_{NM}} (a_{nm}\sqrt{nm} -V_{nm})
\\
&& \quad = \min_{(n,m) \in{\mathcal{K}}_{NM}} T_{nm} \biggl( \frac
{a_{nm}\sqrt{nm} (1+\mathrm{o}(1))}{\sqrt{2(n\log(p^{-1})+m\log(q^{-1}))}} -\sqrt{1+
\frac{\log(NM)}{n\log(p^{-1})+m\log(q^{-1})}} \biggr),
\end{eqnarray*}
which goes to infinity under \nref{assump2} and \nref{cond2a}.
Thus, we have $\gamma_{NM,\mathbf{a}}(\psi_{NM}^{\mathrm{max}}) \to0$.

%
%$$
%$$
%
%Since $N\to\infty,\ n\to\infty, \ N-n\to\infty$ using the Stirling
%formula
%l\in\N,
%we get
%$$
%{N \choose n}\sim\frac{N^{N+1/2}}{n^{n+1/2}(N-n)^{N-n+1/2}\sqrt{2\pi}}.
%$$
%Therefore
%&+&\frac12
%(\log(N/(N-n)-\log(n))+\mathrm{O}(1)\\
%&=&-n\log(p)-(N-n)\log(1-p)-\frac12 (\log(n)+\log(1-p))+\mathrm{O}(1).
%Since $p=n/N\to0$, we have
%$$\log(1-p)\sim-p=\mathrm{o}(1),\ (N-n)\log(1-p)\sim
%-(N-n)p=-n(1-p)=\mathrm{O}(n),
%$$
%and Lemma \ref{L1} follows.

%s6.3 #&#
\subsection{\texorpdfstring{Proof of Proposition \protect\ref{PZ}}{Proof of Proposition 5.1}}\label{PPZ}
It suffices to check that $P_0(\Gamma_{nm}^c)\to0$, where $A^c$ states
for the complement of the event $A$. We have
\begin{eqnarray*}
\Gamma_{nm}^c&=&\bigcup_{C\in\CC_{nm}}
\bigcup_{k_0\le k\le n, l_0\le
l\le m} \bigcup_{V\in\CC_{kl,C}}
\{Y_V>T_{kl}\}\\
&=&\bigcup_{k_0\le
k\le
n, l_0\le l\le m}
\bigcup_{V\in\CC_{kl}}\{Y_V>T_{kl}\}.
\end{eqnarray*}
Since $Y_V \sim\CN(0,1)$ under $P_0$, we have, by definition of
$T_{kl}$ and using the asymptotics $\Phi(-x)\sim
\mathrm{e}^{-x^2/2}/\sqrt{2\uppi}x, x\to\infty$,
\begin{eqnarray*}
P_0\bigl(\Gamma_{nm}^c\bigr)&\le& \sum
_{k_0\le k\le n, l_0\le l\le m} \sum_{V \in\CC_{kl}}
\Phi(-T_{kl})=\sum_{k_0\le k\le n, l_0\le
l\le
m}G_{kl}
\Phi(-T_{kl})
\\
&\le& \sum_{k_0\le k\le n, l_0\le l\le m}\frac{1+\mathrm{o}(1)}{nm
T_{kl}\sqrt{2\uppi}}\to0.
\end{eqnarray*}
Proposition \ref{PZ} follows.

%s6.4 #&#
\subsection{\texorpdfstring{Proof of Proposition \protect\ref{PE}}{Proof of Proposition 5.2}}\label{PPE}
In view of symmetry in $C$, it suffices to check that, for any fixed
$C\in\CC_{nm}$,
\[
E_{0} \biggl(\frac{\mathrm{d}P_{S_C}}{\mathrm{d}P_0}\1_{\Gamma_C}
\biggr)=P_{S_C}(\Gamma_C)\to1,
\]
or, equivalently, $P_{S_C}(\Gamma_C^c)\to0$. Set $z_{kl}^2=a^2kl$. Since
$Y_V\sim\CN(z_{kl},1)$ under $P_{S_C}$ for $V \in\CC_{kl,C}$, we
have
\[
P_{S_C}\bigl(\Gamma_C^c\bigr)\le\sum
_{k_0\le k\le n, l_0\le l\le m} \sum_{V \in\mathcal{C}_{kl,C}}
\Phi(z_{kl}-T_{kl})=\sum_{k_0\le
k\le n, l_0\le
l\le m}G_{kl}^{mn}
\Phi(z_{kl}-T_{kl}),
\]
where $G_{kl}^{mn}=\#(\CC_{kl,C})={n\choose k} {m\choose l}$. Under assumptions
\nref{cond3a} and \nref{cond12oppb} there exists $\delta>0$ such that
%
%e6.2 #&#
\begin{equation}
\label{cond*} b^2=a^2nm<(2-\delta) \bigl(n\log
\bigl(p^{-1}\bigr)+m\log\bigl(q^{-1}\bigr)\bigr).
\end{equation}
Let us show that under \nref{cond*} one has $z_{kl}^2<
T_{kl}^2(1-\delta/2+\mathrm{o}(1))$. In fact, since $\delta_1 n\le k\le n, \delta_1 m\le l\le m$, and by \nref{Tkl}, we have
%
%e6.3 #&#
\begin{eqnarray}\label{ZT}
z_{kl}^2=a^2kl&\le& (2-\delta)
\bigl(k (l/m)\log\bigl(p^{-1}\bigr)+l(k/n)\log\bigl(q^{-1}
\bigr)\bigr)
\nonumber
\\[-8pt]
\\[-8pt]
\nonumber
&\le& (2-\delta) \bigl(k \log\bigl(p^{-1}\bigr)+l\log
\bigl(q^{-1}\bigr)\bigr) \sim (1-\delta/2)T_{kl}^2.
\end{eqnarray}
Thus we get, for some $\delta_2>0$,
\[
\Phi(z_{kl}-T_{kl})\le\exp\bigl(-\delta_2T_{kl}^2
\bigr).
\]
Observe now that, under constraints on $\delta_1 n\le k\le n, \delta_1 m\le l\le m$ we have $\log(G_{kl}^{nm})=\mathrm{O}(n+m)$. This
follows from evaluations similar to the proof of \nref{L1}. On the
other hand, we have $T_{kl}^2\sim2(k\log(p^{-1})+l\log(q^{-1}))\gg
(n+m)$ under the same constraints. This yields
\[
\sum_{k_0\le k\le n, l_0\le l\le
m}G_{kl}^{mn}
\Phi(z_{kl}-T_{kl})\le\sum_{k_0\le k\le n, l_0\le
l\le m}
\exp\bigl(\mathrm{O}(n+m)-\delta_2T_{kl}^2\bigr)\to0.
\]
Proposition \ref{PE} follows. %\hfill\mbox{\ $\Box$}

\subsection{\texorpdfstring{Proof of Lemma \protect\ref{LGen}}{Proof of Lemma 5.1}}\label{LGenP}
The first inequalities in \nref{in1}--\nref{in3} are evident, and we
will prove
the second ones. The proofs are based on the well-known relation: if
$X\sim\CN(0,1)$, then
%
%e6.4 #&#
\begin{equation}
\label{ExpMom} E\bigl(\exp(\tau X)\bigr)=\exp\bigl(\tau^2/2\bigr)\qquad
\forall \tau\in\R.
\end{equation}
Let $V_1=C_1\setminus C_2, V_2=C_2\setminus C_1, V=C_1\cap C_2$,
and observe that the sets $V_1, V_2, V$ are disjoint, $C_1=V_1\cup V, C_2=V_2 \cup V$ and $\#(V_1)=\#(V_2)=nm-kl, \#(V)=kl$.

Let $0<kl<nm$. Let us write the statistics $Y_{C_1}, Y_{C_2}$ in
a more convenient form
\[
Y_{C_1}=\sqrt{1-\rho_{kl}} Y_{V_1}+\sqrt{
\rho_{kl}} Y_V,\qquad Y_{C_2}=\sqrt{1-
\rho_{kl}} Y_{V_2}+\sqrt{\rho_{kl}}
Y_V,
\]
where as above, for $U\subset\{1,\ldots,N\}\times\{1,\ldots,M\}, \#(U)>0$ we set
\[
Y_U=\frac{1}{\sqrt{\#(U)}}\sum_{(i,j)\in U}Y_{ij}.
\]
Observe that $Y_{V_1}, Y_{V_2}, Y_V$ are standard Gaussian and
independent under $P_0$.

Recall that $b=a \sqrt{nm}$ and put $c=b\sqrt{1-\rho_{kl}}$. It is obvious
that $b^2=c^2+z_{kl}^2$. Moreover, by applying \nref{ExpMom}, we get
\begin{eqnarray*}
&& E_0 (\exp\bigl(-b^2+b(Y_{C_1}+Y_{C_2})
\bigr)
\\
&&\quad=E_0 \bigl(\exp\bigl(-c^2/2+cY_{V_1}\bigr)
\bigr)\cdot E_0 \bigl(\exp\bigl(-c^2/2+cY_{V_2}
\bigr) \bigr) \cdot E_0 \bigl(\exp \bigl(-z_{kl}^2+2z_{kl}Y_{V}
\bigr) \bigr)
\\
&&\quad= \exp\bigl(z_{kl}^2\bigr).
\end{eqnarray*}
If $kl=0$ or $kl=nm$, we can prove this in a similar way. Lemma \ref
{LGen} \nref{in1} follows.

In order to get the second inequality, observe that, for $0<kl<nm$
and for any $h\ge0$,
\begin{eqnarray*}
&&E_0 \bigl(\exp\bigl(-b^2+b(Y_{C_1}+Y_{C_2})\bigr)
\1_{Y_{C_1}\le T_{nm}, Y_{C_2}\le T_{nm}} \bigr)
\\
&&\quad\le \mathrm{e}^{-b^2+2 T_{nm}h}E_0 \bigl(\exp \bigl((b-h)
(Y_{C_1}+Y_{C_2})+h(Y_{C_1}+Y_{C_2}-2T_{nm})\bigr)
\1_{Y_{C_1}\le
T_{nm}, Y_{C_2}\le T_{nm}} \bigr)
\\
&&\quad\le \mathrm{e}^{-b^2+ 2T_{nm}h}E_0 \bigl(\exp\bigl((b-h) (Y_{C_1}+Y_{C_2})
\bigr) \bigr)
\\
&&\quad=\mathrm{e}^{-b^2+2T_{nm}h}E_0 \bigl(\exp\bigl((b-h) (1-
\rho_{kl})^{1/2}(Y_{V_1}+Y_{V_2})+2(b-h)
\rho_{kl}^{1/2} Y_V\bigr) \bigr)
\\
&&\quad=\exp\bigl(-b^2+2T_{nm}h+(b-h)^2(1-
\rho_{kl})+2(b-h)^2\rho_{kl}\bigr)
\\
&&\quad=\exp\bigl(-b^2+2T_{nm}h+(b-h)^2(1+
\rho_{kl})\bigr).
\end{eqnarray*}
Taking $h=b-T_{nm}/(1+\rho_{kl})$, we get the second inequality. If
$kl=0$ or $kl=nm$, we can prove this in a similar way. Lemma \ref{LGen}
\nref{in2} follows.

In order to get the third inequality, for
$0<kl<nm$ and for $h\ge0$, we have
\begin{eqnarray*}
&&E_0\bigl (\exp\bigl(-b^2+b(Y_{C_1}+Y_{C_2})\bigr)
\1_{Y_{V}\le T_{kl}} \bigr)
\\
&&\quad= E_0 \bigl(\exp\bigl(-c^2/2 + cY_{V_1}
\bigr) \bigr) \cdot E_0 \bigl(\exp \bigl(-c^2/2 +
cY_{V_2}\bigr) \bigr) \cdot E_0 \bigl(\mathrm{e}^{-z_{kl}^2}
\exp(2z_{kl}Y_{V})\1_{Y_{V}\le
T_{kl}} \bigr)
\\
&&\quad\le \mathrm{e}^{-z_{kl}^2+T_{kl}h}E_0 \bigl(\exp \bigl( (2z_{kl}-h)Y_{V}+h(Y_{V}-T_{kl}
) \bigr)\1_{Y_{V}\le
T_{kl}} \bigr)
\\
&&\quad\le \mathrm{e}^{-z_{kl}^2+T_{kl}h}E_0 \bigl(\exp\bigl((2z_{kl}-h)Y_{V}
\bigr) \bigr)=\exp \bigl(-z_{kl}^2+T_{kl}h+(2z_{kl}-h)^2/2
\bigr).
\end{eqnarray*}
Taking $h=2z_{kl}-T_{kl}$, we get the third inequality. If $kl=nm$,
we can prove this in a similar way. Lemma \ref{LGen} \nref{in3} follows.

%s6.6 #&#
\subsection{\texorpdfstring{Proof of Lemma \protect\ref{L5}}{Proof of Lemma 5.3}}\label{PL5}
Recalling $P_{n,p}(k)={n\choose k} p^k(1-p)^{n-k}$. Using well-known
inequality ${n\choose k}\le(ne/k)^k$, we get
\[
\log\bigl(P_{n,p}(k)\bigr)\le k \bigl(\log(p)+\log(n/k)+1 \bigr).
\]
Since $n/k\le r(p)$, we see that $0\le\log(n/k)\le
\log(r(p))=\mathrm{o}(\log(p^{-1}))$ under the assumption on $r(p)$. This
implies the first relation of Lemma \ref{L5}.

In view of Lemma \ref{L4a}, we have
\[
P_{N,n,n}(k)\le P_{{{\mathcal{HG}}}}(Z\ge k)\le P_{\mathrm{Bin}}(Z\ge k)
\sim P_{n,p}(k).
\]
Lemma \ref{L5} follows. % \hfill\mbox{\ $\Box$}

%s6.7 #&#
\subsection{\texorpdfstring{Proof of Theorem \protect\ref{TUV}}{Proof of Theorem 4.1}}\label{prTUV}
Let us see that $E_0(\hat\sigma^2) = \sigma^2$ and that
$\Var_0(\hat\sigma^2) = 2\sigma^4/(NM)$. Denote by
\[
O_B = \biggl\{\biggl\llvert \frac{\hat\sigma^2}{\sigma^2} - 1\biggr\rrvert \leq
\frac B{\sqrt{NM}} \biggr\},
\]
with $B \to\infty$ such that $B/\sqrt{NM} \to0$. Then,
$P_0((O_B)^c)\leq2B^{-2} = \mathrm{o}(1)$.

It is easy to see that $P_0(\hat t_{\mathrm{lin}}>H) \leq P_0(t_{\mathrm{lin}}/\sigma
>H \hat\sigma/\sigma, O_B)+ P_0((O_B)^c) = P_0(t_{\mathrm{lin}}/\sigma>
\tilde H )+ \mathrm{o}(1)=\mathrm{o}(1)$, as $\tilde H^2 = H^2 (1-B/\sqrt{NM} )\to
\infty$.

Similarly, for $\hat t_{\mathrm{max}}$, we put $T_{nm,\delta} =
T_{nm}\sqrt{1+\delta/2} = $ and $\tilde T_{nm}^2 = T^2_{nm,\delta}
(1-B\sqrt{NM})$ and then $P_0(\hat t_{\mathrm{max}} >T_{nm,\delta}) \leq
P_0(t_{\mathrm{max}} /\sigma> \tilde T_{nm}) + P_0((O_B)^c) =\mathrm{o}(1)$, for our
choice of $T_{nm,\delta}$. These imply $\a(\hat\psi^*)\to0$.

Under the alternative, let us see that
\[
\b\bigl(\hat\psi^*, S\bigr) \leq\min\bigl\{ P_S(\hat
t_{\mathrm{lin}} \leq H), P_S(\hat t_{\mathrm{max}} \leq
T_{nm,\delta}) \bigr\},
\]
and it suffices to check that either $P_S(\hat t_{\mathrm{lin}} \leq H)\to0$
or $P_S(\hat t_{\mathrm{max}} \leq T_{nm,\delta})\to0$. We can decompose
\begin{eqnarray*}
\hat\sigma^2 &=& \frac{\sigma^2}{NM} \biggl( \sum
_{(i,j)\notin
C} \xi_{ij}^2 + \sum
_{(i,j)\in C} (s_{ij}+\xi_{ij})^2
\biggr)
\\
&=& \frac{\sigma^2}{NM} \biggl(\sum_{(i,j)}
\xi_{ij}^2 + 2 \sum_{(i,j)\in
C}
\frac{s_{ij}}\sigma\xi_{ij}+ \sum_{(i,j)\in C}
\frac
{s_{ij}^2}{\sigma^2} \biggr).
\end{eqnarray*}
We get
\[
E_S\bigl(\hat\sigma^2\bigr) = \sigma^2
(1+G), \qquad \Var_S\bigl(\hat\sigma^2\bigr) =
\frac{2\sigma^4}{NM}(1+ 2G)\qquad\mbox{where } G = \frac1{\sigma^2 NM}\sum
_{(i,j) \in C} s_{ij}^2.
\]
Denote by $R:=B \sqrt{\frac2{NM} (1+2G)}$ with $B \to\infty$ such
that $B = \mathrm{o}( \sqrt{NM})$, and by
\[
O_{SB} = \biggl\{\biggl\llvert \frac{\hat\sigma^2}{\sigma^2} - 1 -G\biggr\rrvert
\leq R \biggr\}.
\]
Then $P_S((O_{SB})^c)\leq B^{-2} = \mathrm{o}(1)$ and $R=\mathrm{o}(1+G)$. Recalling
\[
E_S(t_{\mathrm{lin}}/\sigma)=(NM)^{-1/2}\sum
_{(i,j) \in C} s_{ij}/\sigma\ge a\sqrt{nmpq}/\sigma,
\]
we see that
\[
G\le\frac{\max_{(i,j) \in C}s_{ij}}{\sum_{(i,j) \in C}s_{ij}} \frac{  (\sum_{(i,j) \in C} s_{ij} )^2}{\sigma^2NM}=\tau_{G} \bigl(
E_S(t_{\mathrm{lin}}/\sigma) \bigr)^2.
\]
This implies
%
%e6.5 #&#
\begin{equation}
\label{GR} 1+G+R\le\bigl(1+\tau_{G}\bigl(E_S(t_{\mathrm{lin}}/
\sigma)\bigr)^2\bigr) \bigl(1+\mathrm{o}(1)\bigr).
\end{equation}

Let $E_S(t_{\mathrm{lin}}/\sigma)\to\infty$, which holds when
$a\sqrt{nmpq}/\sigma\to\infty$. By our choice of $H$, \nref{altmax}
and~\nref{GR} we have $H \sqrt{1+G+R}=\mathrm{o}(E_S(t_{\mathrm{lin}}/\sigma)).$
Applying the Chebyshev inequality and since
$\Var_S(t_{\mathrm{lin}}/\sigma)=1$, we have
\begin{eqnarray*}
P_S(\hat t_{\mathrm{lin}} \leq H) &\leq& P_S(t_{\mathrm{lin}}/
\sigma\leq H \hat\sigma /\sigma, O_{SB}) + P_S
\bigl((O_{SB})^c\bigr)
\\
&\leq& P_S(t_{\mathrm{lin}}/\sigma\leq H \sqrt{1+G+R} ) +\mathrm{o}(1)
\\
&\leq& P_S\bigl(E_S(t_{\mathrm{lin}}/\sigma)
-t_{\mathrm{lin}}/\sigma\leq E_S(t_{\mathrm{lin}}/\sigma )- H
\sqrt{1+G+R}\bigr) +\mathrm{o}(1)
\\
& \leq& \bigl(E_S(t_{\mathrm{lin}}/\sigma) - H \sqrt{1+G+R}
\bigr)^{-2} +\mathrm{o}(1)=\mathrm{o}(1).
\end{eqnarray*}
This proves that $P_S(\hat t_{\mathrm{lin}} \leq H)
\to0$ if $ E_S(t_{\mathrm{lin}}/\sigma)\to\infty$.

{\spaceskip=0.18em plus 0.05em minus 0.02em If $ E_S(t_{\mathrm{lin}}/\sigma) = \mathrm{O}(1)$ (this is possible when
$a\sqrt{nmpq}/\sigma=\mathrm{O}(1)$ only), we have $\sqrt{1+G+R } =$} $1+\mathrm{o}(1)$
by \nref{altmax} and \nref{GR}. Therefore
\begin{eqnarray*}
P_S(\hat t_{\mathrm{max}} \leq T_{nm,\delta}) &
\leq&P_S(t_{\mathrm{max}} \leq T_{nm,\delta} \hat\sigma/\sigma,
O_{SB})+ P_S\bigl((O_{SB})^c\bigr)
\\
&\leq& P_S(t_{\mathrm{max}} \leq T_{nm,\delta} \sqrt{1+G+R})
+\mathrm{o}(1)=\mathrm{o}(1)
\end{eqnarray*}
for our choice of $T_{nm,\delta}$ and by the second assumption
\nref{42} (compare with the proof of Theorem~\ref{TU}). These
implies $\b_{nm,a}(\hat\psi^*)\to0$. Thus $\g_{nm,a}(\hat\psi^*)\to
0$.

%s6.8 #&#
\subsection{\texorpdfstring{Proof of the upper bound of Theorem~\protect\ref{TE}}
{Proof of the upper bound of Theorem 4.2}}\label{sectionTE}
It follows the same lines as that of Theorem~\ref{TU}.
We use Markov inequality and
bound from above exponential moments of our test statistics
(as they are not having Gaussian distribution in this case).

We use repeatedly the well-known facts that, $A'(s^0)=0$ and $A''(s^0)
= 1$ for
centered and reduced random variable at $s^0$. Moreover,
\[
E_s\bigl[\mathrm{e}^{u Y}\bigr] = \mathrm{e}^{A(s+u)-A(s)}
\]
for any $s$ and $u$ such that $s$ ans $s+u$ are interior points of the
parameter space.
For the statistic $t_{\mathrm{lin}}$, we have
\begin{eqnarray*}
\alpha\bigl(\psi_H^{\mathrm{lin}}\bigr) &= & P_{s^0}
\bigl(t_{\mathrm{lin}}H >H^2\bigr) \leq \mathrm{e}^{-H^2}
E_{s^0}\bigl[\mathrm{e}^{t_{\mathrm{lin}} H}\bigr]
\\
& \leq& \mathrm{e}^{-H^2} \prod_{(i,j)} E_{s^0}
\bigl[\mathrm{e}^{Y_{ij} H/\sqrt{NM}}\bigr]
\\
&\leq& \exp\biggl( \biggl(A\biggl(s^0+\frac H{\sqrt{NM}}\biggr)-A
\bigl(s^0\bigr)\biggr)NM - H^2\biggr).
\end{eqnarray*}
For $H \to\infty$ as in our theorem $H/\sqrt{NM} \to0$, then we get
$A(s^0+\frac H{\sqrt{NM}}) - A(s^0)=\frac{H^2}{2 NM}(1+\mathrm{o}(1))$ and
$\alpha(\psi_H^{\mathrm{lin}}) \leq c\mathrm{e}^{-H^2/2(1-\mathrm{o}(1))} \to0$, for some
constant $c>0$.

Under the alternative,
\begin{eqnarray*}
\beta\bigl(\psi_H^{\mathrm{lin}},S\bigr) &=& P_S[-t_{\mathrm{lin}}
+H \geq0] \leq \mathrm{e}^H E_S\bigl[\mathrm{e}^{-t_{\mathrm{lin}}}\bigr]
\\
& \leq& \mathrm{e}^H \prod_{(i,j)\notin C} E_{s^0}
\bigl[\mathrm{e}^{-Y_{ij} /\sqrt{NM}}\bigr] \prod_{(i,j)\in C} E_{s_{ij}}
\bigl[\mathrm{e}^{-Y_{ij} /\sqrt{NM}}\bigr]
\\
& \leq& \mathrm{e}^{ H + (NM-nm)  (A(s^0 -
(1/{\sqrt{NM}}))-A(s^0) )}\prod_{(i,j)\in C} \mathrm{e}^{ A(s_{ij} -
(1/{\sqrt{NM}}) ) - A(s_{ij})}.
\end{eqnarray*}
On the one hand, $A(s^0 - \frac1{\sqrt{NM}})-A(s^0) \sim-\frac1{2NM}$.
On the other hand, it is easy to check that $A''(s) \geq0$.
Thus, $A$ is a convex function and
$A'$ is increasing. This implies
$ A(s_{ij} - \frac1{\sqrt{NM}}) - A(s_{ij}) \leq-\frac1{\sqrt{NM}}
A'(s_{ij})$ which is
less or equal than $-\frac1{\sqrt{NM}} A'(s^0+a)$ under the
alternative. Finally,\vspace*{-2pt}
\begin{eqnarray*}
\beta\bigl(\psi_H^{\mathrm{lin}},S\bigr) &\leq& \exp\biggl(H -
\frac12(1-pq) - A'\bigl(s^0+a\bigr) \frac
{nm}{\sqrt{NM}}
\biggr)
\\[-3pt]
&\leq& c \exp\bigl(H- A'\bigl(s^0+a\bigr) \sqrt{nmpq}
\bigr) \to0\vspace*{-2pt}
\end{eqnarray*}
under our assumption.

For the statistic $t_{\mathrm{max}}$,\vspace*{-1pt}
\begin{eqnarray*}
\alpha\bigl(\psi^{\mathrm{max}}\bigr) &= & P_{s^0}(t_{\mathrm{max}} >
T_{nm,\delta}) \leq\sum_C P_{s^0}
\bigl(t_{\mathrm{max}} T_{nm,\delta} > T^2_{nm,\delta}\bigr)
\\[-2pt]
& \leq& \sum_C \mathrm{e}^{-T_{nm,\delta}^2}E_{s^0}
\bigl[\mathrm{e}^{Y_{C} T_{nm,\delta}}\bigr] \leq\sum_C
\mathrm{e}^{-T_{nm,\delta}^2}E_{s^0}\bigl[\mathrm{e}^{\sum_{(i,j) \in C}Y_{ij}
{T_{nm,\delta}}/{\sqrt{nm}}}\bigr]
\\[-2pt]
&\leq& G_{nm} \mathrm{e}^{-T_{nm,\delta}^2} \mathrm{e}^{ ( A(s^0+
{T_{nm,\delta
}}/{\sqrt{nm}}) -A(s^0) ) \times nm}.
\end{eqnarray*}
As $T_{nm,\delta}/\sqrt{nm}\asymp ( \log(p^{-1})/m + \log
(q^{-1})/n )^{1/2}\to0$,
we have $A(s^0+T_{nm,\delta}/\sqrt{nm}) -A(s^0) \asymp T^2_{nm,\delta}/(2nm)$
and this gives\vspace*{-1pt}
\[
\alpha\bigl(\psi^{\mathrm{max}}\bigr) \leq G_{nm} \mathrm{e}^{-({T_{nm,\delta}^2}/2) (1-\mathrm{o}(1))}
\to0.\vspace*{-1pt}
\]
For the second-type error,\vspace*{-2pt}
\begin{eqnarray*}
\beta\bigl(\psi^{\mathrm{max}},S\bigr) &\leq& P_S[-
Y_C + T_{nm,\delta} \geq0] \leq \mathrm{e}^{T_{nm,\delta}}
E_S\bigl[\mathrm{e}^{-Y_C}\bigr]
\\[-3pt]
& \leq& \exp(T_{nm,\delta}) \cdot\prod_{(i,j) \in C} E_{s_{ij}}
\bigl[\mathrm{e}^{-Y_{ij} /\sqrt{nm}}\bigr]
\\[-3pt]
& \leq& \exp \biggl( T_{nm,\delta} + \sum_{(i,j)\in C}
\biggl[A\biggl(s_{ij}- \frac 1{\sqrt{nm}}\biggr)-A(s_{ij})
\biggr] \biggr)
\\[-3pt]
&\leq& \exp \biggl(T_{nm,\delta} -nm A'\bigl(s^0
+a\bigr) \frac1{\sqrt {nm}} \biggr) \to0,
\end{eqnarray*}
by the choice of $\delta> 0$ small enough.\vspace*{-2pt}

%s6.9 #&#
\subsection{\texorpdfstring{Proof of Theorem~\protect\ref{TTS}}{Proof of Theorem 4.3}}\label{ProofTS}

%s6.9.1 #&#
\subsubsection{Proof of the upper bounds}

We have, under the null hypothesis, $z_{\mathrm{lin}} \sqrt{2NM}$ has a $\chi^2$
distribution with $E_0(z_{\mathrm{lin}}) = 0$ and
$\operatorname{Var}_0(z_{\mathrm{lin}}) = 1$.
This implies that $P_0(z_{\mathrm{lin}} >H) \to0$ as $H \to\infty$.

For $z_{\mathrm{max}}$ we will use the moment generating function
of the $\chi^2$ distribution. We have\vspace*{-2pt}
\begin{eqnarray*}
P_0(z_{\mathrm{max}} >T_{nm,\delta}) &\leq& \sum
_C P_0(Z_C >T_{nm,\delta} )
\leq G_{nm} P_0\bigl(T_{nm,\delta} Z_C
>T^2_{nm,\delta}\bigr)
\\[-2pt]
&\leq& G_{nm} \mathrm{e}^{-T^2_{nm,\delta}} E_0\bigl(\mathrm{e}^{ T_{nm,\delta} Z_C}
\bigr)
\\
%&\leq& G_{nm} \mathrm{e}^{-T_{nm,\delta}^2- T_{nm,\delta} \sqrt{\frac{nm}2} }
%E_0\left(\exp(\frac{T_{nm,\delta}}{\sqrt{2nm}} \sum_C Y_{ij}^2)\right)\
&\leq& G_{nm} \mathrm{e}^{-T_{nm,\delta}^2- T_{nm,\delta} \sqrt{{nm}/2} }
\biggl(1-\frac{2 T_{nm,\delta}}{\sqrt{2nm}}\biggr)^{-nm/2}
\\
&\leq& G_{nm} \mathrm{e}^{-({T_{nm,\delta}^2}/2) (1+\mathrm{o}(1))} = \mathrm{o}(1),
\end{eqnarray*}
by the choice of $T_{nm,\delta}$ in our theorem. Indeed, $2
T_{nm,\delta
}^2/(nm) = \mathrm{O}(a^4) \to0$, by assumption~\nref{tend}.

Again, $\beta(\psi^z, S) \leq\min\{P_S(z_{\mathrm{lin}} \leq H), P_S(z_{\mathrm{max}}
\leq T_{nm,\delta})\}$.
Under the alternative, $S=S_C$ and $z_{\mathrm{lin}}$ has mean
$E_S(z_{\mathrm{lin}})=\lambda/\sqrt{NM}$ and variance
$\operatorname{Var}_S(z_{\mathrm{lin}})=1+2\lambda/ (NM)$, where $\lambda= \sum_C s_{ij}^2 $.
We have,
\[
\frac{\lambda}{\sqrt{2NM}} \geq\frac{a^2 nm}{\sqrt{2NM}} \geq \frac
{a^2}{\sqrt{2}} \sqrt{nmpq}.
\]
Therefore, if $a^2\sqrt{nmpq} \to\infty$, we have
\begin{eqnarray*}
P_S(z_{\mathrm{lin}} \leq H) &\leq&\frac{\operatorname{Var}_S(z_{\mathrm{lin}})}{(E_S(z_{\mathrm{lin}})-H)^2} =
\frac{1+2\lambda/(NM)}{(1-c)^2 \lambda^2 /(2NM)}
\\
& \leq& \frac2{(1-c)^2 a^4 nmpq} + \frac{4}{(1-c)^2 a^2 nm} =
\mathrm{o}(1).
\end{eqnarray*}

Under the alternative,
\begin{eqnarray*}
P_S(z_{\mathrm{max}} \leq T_{nm,\delta}) &\leq&
\mathrm{e}^{T_{nm,\delta}} E_{S_C}\bigl[\mathrm{e}^{-Z_C}\bigr] \leq
\mathrm{e}^{T_{nm,\delta} +\sqrt{nm/2}} E_{S_C} \biggl[\exp \biggl( t \sum
_C Y_{ij}^2 \biggr)\biggr],
\end{eqnarray*}
where $t = -1/\sqrt{2nm} < 1/2$. Therefore,
\begin{eqnarray*}
E_{S_C} \biggl[\exp \biggl( t \sum_C
Y_{ij}^2 \biggr) \biggr] & = & \exp \biggl(
\frac{\lambda t}{1-2t} -\frac{nm}2 \log(1 -2t) \biggr)
\\
&\leq& \exp \biggl( -\frac{\lambda}{\sqrt{2nm}} \frac 1{1+\sqrt{2/(nm)}} -
\frac{nm}2 \biggl(\sqrt{\frac2{nm}} - \frac1{nm}\biggr) \biggr)
\\
& \leq& \exp \biggl(- \frac{a^2 \sqrt{nm}}{\sqrt{2}} \frac 1{1+\sqrt{2/(nm)}} -\sqrt{
\frac{nm}2 } +\frac12 \biggr).
\end{eqnarray*}
In conclusion, if $\liminf a^4 nm /(4(n\log(p^{-1}) + m\log(q^{-1})))
>1$ we have
\[
P_{S_C}(z_{\mathrm{max}}\leq T_{nm,\delta}) \leq\sqrt{e} \exp
\biggl(T_{nm,\delta
} -\frac{a^2 \sqrt{nm}}{\sqrt{2}} \frac1{1+\sqrt{2/(nm)}} \biggr) \to0.
\]

%s6.9.2 #&#
\subsubsection{Proof of the lower bounds}

We follow the lines of the proof of Theorem~\ref{TL}. The prior on the
set of matrices is $\pi= G_{nm}^{-1} \sum_{C \in\mathcal{C}_{nm}}
\pi_C$, where, under $\pi_C$, the matrix $S=S_C$ has $s_{ij}=0$ with
probability 1 for all $(i,j) \notin C$ and $s_{ij}$ is either $a$ and
$-a$ with probability $1/2$, for all $(i,j) \in C$.

Let $P_{S_C}$ denote the likelihood of the random variables in $\mathbf
{Y}$ when $S=S_C$ and $P_\pi$ denote the mixture of likelihoods $P_\pi
= G_{nm}^{-1} \sum_{C \in\mathcal{C}_{nm}}P_{S_C}$. Therefore, the
likelihood ratio $L_\pi(Y)$ is
\begin{eqnarray*}
L_\pi(Y) & = & \frac{\mathrm{d}P_\pi}{\mathrm{d} P_0}(Y) = \frac1{G_{nm}} \sum
_{C \in\mathcal{C}_{nm}} \frac{\mathrm{d}P_{\pi_C}}{\mathrm{d}
P_0}(Y)
\\
& = & \frac1{G_{nm}} \sum_{C \in\mathcal{C}_{nm}}
\prod_{(i,j) \in C} \mathrm{e}^{-a^2/2} \cosh(a Y_{ij})
\\
& = & \frac1{G_{nm}} \sum_{C \in\mathcal{C}_{nm}}
\mathrm{e}^{-a^2 nm/2} \exp \biggl( \sum_{(i,j) \in C} \log\bigl(
\cosh( a Y_{ij})\bigr) \biggr) .
\end{eqnarray*}
Note that $E_0[\cosh(aY_{ij})] = \mathrm{e}^{a^2/2}$ and that
$E_0[\cosh^2(aY_{ij})] = (1+\mathrm{e}^{2a^2})/2$.
% = 1+ a^2 + \mathrm{o}(a^2)$ as $a\to0$. This implies that w
We can reproduce the proof of
Theorem~\ref{TL} with $a^2$ replaced by $a^4/2$. For example, in
\nref{in1} we have
\[
g(k,l) = E_0\biggl[\exp\biggl(-a^2nm + \sum
_{C_1} \log\bigl(\cosh(a Y_{ij})\bigr) +\sum
_{C_2} \log\bigl(\cosh(a Y_{ij})\bigr)
\biggr) \1_{\Gamma_{C_1} \cap\Gamma_{C_2}}\biggr].
\]
We can show, as in the proof of \nref{in1}, that
\begin{eqnarray*}
g(k,l) &\leq& \mathrm{e}^{-a^2 nm } E_0\biggl[\exp\biggl( \sum
_{V_1 \cup V_2} \log\bigl(\cosh (a Y_{ij})\bigr) +2 \sum
_{V} \log\bigl(\cosh(a Y_{ij})\bigr)
\biggr)\biggr]
\\
&\leq& \mathrm{e}^{-a^2 nm } \mathrm{e}^{2 ({a^2}/2) (nm-kl)} E_0^{kl}\bigl[
\cosh^2(a Y)\bigr] \leq \mathrm{e}^{-a^2 kl} \biggl( \frac{1+\mathrm{e}^{2a^2}}2
\biggr)^{kl}
\\
&=&\bigl(\cosh\bigl(a^2\bigr)\bigr)^{kl} \leq
\mathrm{e}^{({a^4}/2) kl},
\end{eqnarray*}
where $V_1,  V_2$ and $V$ are defined in the proof of Lemma~\ref{LGen}.

The relations \nref{in2} and \nref{in3} could be replaced by the
following:
%
%e6.6 #&#
%e6.7 #&#
\begin{eqnarray}
\label{in2a} g(k,l) &\le& \exp \biggl(-(T_{nm}-b)^2+
\frac{\rho_{kl}T_{nm}^2}{1+\rho_{kl}}+\mathrm{o}\bigl(T_{nm}^2\bigr) \biggr),
\\
\label{in3a} g(k,l)&\le&\exp \bigl(T^2_{kl}/2-(T_{kl}-z_{kl})^2+\mathrm{o}
\bigl(T_{kl}^2\bigr) \bigr),
\end{eqnarray}
where $b^2=nm a^4/2, z_{kl}^2=b^2\rho_{kl}$, under the same
constraints. The inspection of the proofs of~\nref{in2} and
\nref{in3} shows that, in order to prove \nref{in2a} and
\nref{in3a}, one could use the following relation in place of
\nref{ExpMom}:
%
%e6.8 #&#
\begin{equation}
\label{cosh} E_0\bigl[\mathrm{e}^{\tau(\log(\cosh(aY)) - {a^2}/2 +{a^4}/4)}\bigr] = \exp \biggl(
\frac{\tau^2 a^4}{4}+ \mathrm{o}\bigl(a^4\bigr) \biggr)
\end{equation}
for $a\to0$ and $\tau\in\mathbb{R}^+, \tau=\mathrm{O}(1)$.

In order to prove \nref{cosh}, we can split the expected value over
the events $\{\tau a^2Y^2
>\delta^2\}$ and $\{\tau a^2Y^2 \leq\delta^2\}$, respectively, for
some small
enough $\delta>0$ such that $\delta/a\sqrt{\tau} \to\infty$ (we
choose $ \delta= (\tau a^2 )^{1/4}$). Firstly, we use the
inequality $\cosh(x) \leq \mathrm{e}^{x^2/2}$ and get
\begin{eqnarray*}
E_0\bigl[\mathrm{e}^{\tau\log(\cosh(aY))} \cdot\1_{\tau a^2 Y^2 > \delta^2}\bigr] & \leq&
E_0\bigl[\mathrm{e}^{\tau a^2 Y^2/2} \cdot\1_{\tau a^2 Y^2 > \delta^2}\bigr]
\\
& \leq& 2 \int_{\delta/a\sqrt{\tau}}^\infty \mathrm{e}^{-(1-\tau a^2) y^2
/2}
\frac{\mathrm{d}y}{\sqrt{2\uppi}}
\\
&\leq& \sqrt{\frac2\uppi} \frac{\delta/(a\sqrt{\tau})}{1 - \tau
a^2} \exp\biggl(-\frac{1-\tau a^2}{2}
\frac{\delta^2}{\tau a^2}\biggr) =\mathrm{o}\bigl(\tau^2 a^4\bigr).
\end{eqnarray*}
Secondly, on the event $\{\tau a^2Y^2 \leq\delta^2\}$ we use the
Taylor expansions $\log(\cosh(x))=x^2/2-x^4/12(1+\mathrm{o}(1))$, $\mathrm{e}^x=1+x+x^2/2(1+\mathrm{o}(1))$,
$x=\mathrm{o}(1)$. Denote $ U = \log(\cosh(aY)) -
E_0(\log(\cosh(aY))). $
We have
%
%e6.9 #&#
\begin{equation}
\label{cosh1} E_0\bigl(\log\bigl(\cosh(aY)\bigr)\bigr) =
\frac{a^2}2 - \frac{a^4}4\bigl(1+\mathrm{o}(1)\bigr), \qquad\Var_0(U) =
\frac{a^4}2\bigl(1+\mathrm{o}(1)\bigr),
\end{equation}
and, since $\tau U=\mathrm{o}(1)$,
\begin{eqnarray*}
E_0 \bigl[\mathrm{e}^{\tau U} \cdot\1_{\tau a^2Y^2\leq\delta^2} \bigr]
% &=& E_0\left[\mathrm{e}^{\tau U} \cdot\1_{\tau U\leq\frac{\delta^2}2 - \frac{
&= & E_0 \biggl[\biggl(1+ \tau U +
\frac{\tau^2 U^2}2 \bigl(1+\mathrm{o}(1)\bigr)\biggr)\cdot\1_{\tau
a^2Y^2\leq\delta^2} \biggr]
\\
&=& 1+ \frac{\tau^2 \Var_0(U)}{2} \bigl(1+\mathrm{o}(1)\bigr)
\\
&&{} -  E_0 \biggl[\biggl(1+ \tau U + \frac{\tau^2 U^2}2 \bigl(1+\mathrm{o}(1)
\bigr)\biggr) \cdot\1_{\tau a^2Y^2> \delta^2} \biggr].
\end{eqnarray*}
The last expected value is $\mathrm{o}(\tau^2 a^4)$ and this gives
\begin{eqnarray*}
E_0\bigl[\mathrm{e}^{\tau\log(\cosh(aY))} \cdot\1_{\tau a^2 Y^2 \leq
\delta^2}\bigr] & = &
\mathrm{e}^{\tau E_0[\log\cosh(a Y)]+ ({\tau^2 a^4}/4)
(1+\mathrm{o}(1))}. % &\leq& \mathrm{O}(1) \exp(\frac{\tau a^2}2 - \frac{\tau a^4}4 +\frac{\tau^2
%a^4}4 ).
\end{eqnarray*}
Together with the first relation in \nref{cosh1}, this ends the proof of
\nref{cosh}.

%s6.10 #&#
\subsection{\texorpdfstring{Proof of Theorem \protect\ref{T2}}{Proof of Theorem 4.4}}\label{PT2}

%s6.10.1 #&#
\subsubsection{Proof of the lower bounds}\label{PT2L}

Let $K=[N/n],  L=[M/m]$ and consider only non-overlapping rectangles
\[
R_{kl}=\bigl\{(i,j)\dvt  n(k-1)+1\le i\le nk, m(l-1)+1\le j\le m l\bigr\},\qquad
1\le k\le K, 1\le l\le L.
\]
Let $S_{kl}$ be the matrix with the
elements $s_{ij}=0$ if $(i,j)\notin R_{kl}$ and $s_{ij}=a$ if
$(i,j)\in R_{kl}$. Consider the prior
\[
\pi=\frac{1}{KL}\sum_{k=1}^K\sum
_{l=1}^L\delta_{S_{kl}}.
\]
By construction, $\pi(\{S_{kl},  k, l \})=1$. The likelihood ratio is
of the
form
\[
L(Y)=\frac{\mathrm{d}P_\pi}{\mathrm{d}P_0}(Y)=\frac{1}{KL}\sum_{k=1}^K
\sum_{l=1}^L\frac
{\mathrm{d}P_{S_{kl}}}{\mathrm{d}P_0}(Y)=
\frac{1}{KL}\sum_{k=1}^K\sum
_{l=1}^L\exp\bigl(-b^2/2+bZ_{kl}
\bigr),
\]
where
\[
Z_{kl}=\frac{1}{\sqrt{nm}}\sum_{(i,j)\in R_{kl}}Y_{ij},\qquad
b^2=nma^2.
\]
Note that $Z_{kl}\sim\CN(0,1)$ under $P_0$ and are independent in
$k,l$. It is sufficient to check that
$L(Y)\to1$ in $P_0$-probability. Let us consider the truncated
likelihood ratio
\[
\tilde L(Y)=\frac{1}{KL}\sum_{k=1}^K
\sum_{l=1}^L\exp\bigl(-b^2/2+bZ_{kl}
\bigr)\1_{Z_{kl}<T_{KL}},
\]
where we set
\[
T_{KL}=\sqrt{2\log(KL)}\sim\sqrt{2\bigl(\log\bigl(p^{-1}
\bigr)+\log\bigl(q^{-1}\bigr)\bigr)}.
\]
Since
\[
P_0(L \ne\tilde L)\le\sum_{k=1}^K
\sum_{l=1}^L P_0(Z_{kl}
\ge T_{KL})\to0,
\]
it suffices to check that $\tilde L(Y)\to1$ in $P_0$-probability.

Observe now that $T_{KL}-b\to\infty$ under the assumptions of
the theorem, and it suffices to consider the case $b>c T_{kl}$ for some
$c\in(1/2,1)$. We have
%(compare with the calculations in the proof
%of Lemma \ref{L2bis} in Section \ref{PL2bis})
%
\begin{eqnarray*}
E_0\bigl(\tilde L(Y)\bigr)&=&\frac{1}{KL}\sum
_{k=1}^K\sum_{l=1}^L
E_0\bigl(\exp\bigl(-b^2/2+bZ_{kl}\bigr)
\1_{Z_{kl}<T_{KL}}\bigr)=\Phi(T_{KL}-b)\to1,
\\
\Var_0\bigl(\tilde L(Y)\bigr)&=&\frac{1}{(KL)^2}\sum
_{k=1}^K\sum_{l=1}^L
\Var_0\bigl(\exp \bigl(-b^2/2+bZ_{kl}\bigr)
\1_{Z_{kl}<T_{KL}}\bigr)
\\
&\le& \frac{1}{(KL)^2}\sum_{k=1}^K
\sum_{l=1}^L E_0 \bigl(\exp
\bigl(-b^2+2bZ_{kl}\bigr)\1_{Z_{kl}<T_{KL}}\bigr)
\\
&=&\frac{1}{KL}\exp\bigl(b^2\bigr)\Phi(T_{KL}-2b)
\\
&\le& \exp\bigl(b^2-(T_{KL}-2b)^2/2-T_{KL}^2/2
\bigr)
\\
&=&\exp\bigl(-(T_{KL}-b)^2\bigr)\to0.
\end{eqnarray*}
Proof of the lower bounds in Theorem \ref{T2}.

%s6.10.2 #&#
\subsubsection{Proof of the upper bounds}\label{PT2U}

Set $T_{KL}=\sqrt{2\log(KL)}$ and observe that, by the choice of
$\eta$ and since $pq\to0$, we have
\begin{eqnarray*}
T_{KL}&=&\sqrt{2\log(NM/nm)-4\log(\eta)+\mathrm{o}(1)}
\\
&=& \bigl(2\log \bigl((pq)^{-1} \bigr) \bigr)^{1/2} \bigl(1+
\bigl(\log (\eta )+\mathrm{o}(1) \bigr)/\log \bigl((pq)^{-1} \bigr)
\bigr)^{1/2}
\\
&\sim& \sqrt{2\bigl(\log\bigl(p^{-1}\bigr)+\log\bigl(q^{-1}
\bigr)\bigr)}.
\end{eqnarray*}

For type I errors, we have
\[
\a(\psi_Z)\le\sum_{k=1}^K
\sum_{l=1}^L P_0(Z_{kl}>T_{KL})=KL
\Phi(-T_{KL})\to0.
\]

Let the alternative $S_E$ correspond to the matrix with entry $a>0$ at
positions in $E=E_{k^*l^*}$
and~0 elsewhere. As previously, $E=E_{k^*l^*}, 0\le
k^*\le N-n, 0\le l^*\le M-m$ consists of $(i,j)$ such that $
k^*<i\le k^*+n, l^*<i\le l^*+m$. By construction, we can take
$k,l, 1\le k\le K, 1\le l\le L$ such that $|n_k-k^*|\le n\eta, |m_l-l^*|\le m\eta$.
Therefore, the matrix $E_{k^*l^*}$ will overlap
with the matrix $E_{n_{k} m_{l}}$ from our test procedure significantly:
\begin{eqnarray*}
\tilde n&=&\# \bigl(\bigl\{k^*+1,\ldots,k^*+n\bigr\}\cap\{n_k+1,\ldots,n_k+n
\} \bigr)\ge n(1-\eta),
\\
\tilde m&=&\# \bigl(\bigl\{l^*+1,\ldots,l^*+m\bigr\}\cap\{m_l+1,\ldots,m_l+m
\} \bigr)\ge m(1-\eta).
\end{eqnarray*}
Observe that
\[
\b(\psi_Z, S_E)\le P_{S_E}(Z_{kl}
\le T_{KL}).
\]
Moreover, $Z_{kl}\sim\CN(\tilde b, 1)$ under $P_{S_E}$
where we recall that $b=a\sqrt{nm}$ and we put
\[
\tilde b=\frac{a \tilde n \tilde m}{\sqrt{nm}}\ge b(1-\eta)^2\sim b.
\]
This yields
\[
\b(\psi_Z, S_E)\le\Phi\bigl(T_{KL}-b
\bigl(1+\mathrm{o}(1)\bigr)\bigr) \to0
\]
under assumptions of theorem. Proof of the upper bounds in Theorem \ref{T2} follows.
\end{appendix}

\section*{Acknowledgements}
Research was
partially supported by RFBR Grant 11-01-00577 and by the Grant
NSh--4472.2010.1. The author acknowledges support from the CNRS for his visit to the University
Paris-Est Marne-la-Vall\'ee.

% imsref loaded by akundreckaite, 2012-09-21 11:03:15
%

%suskaldyti doi

\printhistory

\end{document}